\documentclass[11pt]{article}
\usepackage[T1]{fontenc}
\usepackage[utf8]{inputenc}
\usepackage{lmodern}
\usepackage[margin=1in]{geometry}
\usepackage{amsmath,amssymb,amsthm,mathtools}
\usepackage{microtype,booktabs,array,xcolor,graphicx}
\usepackage[numbers,sort&compress]{natbib}
\usepackage[colorlinks=true,linkcolor=blue!45!black,citecolor=blue!45!black,urlcolor=blue!45!black]{hyperref}
\hypersetup{pdftitle={The Marvelous Meixner Polynomials: Second Kind, First Class},pdfauthor={Olga V. Holtz},pdfsubject={Meixner-Pollaczek polynomials, Carlitz, Wall fractions, spectral duality, and Cauchy couplings}}
\providecommand{\doi}[1]{\href{https://doi.org/#1}{doi: \nolinkurl{#1}}}
\newtheorem{theorem}{Theorem}[section]
\newtheorem{proposition}[theorem]{Proposition}
\newtheorem{corollary}[theorem]{Corollary}
\theoremstyle{remark}
\newtheorem{remark}[theorem]{Remark}
\newcommand{\R}{\mathbb R}
\newcommand{\C}{\mathbb C}
\newcommand{\E}{\mathbb E}
\newcommand{\Pbb}{\mathbb P}
\newcommand{\dd}{\,\mathrm d}
\numberwithin{equation}{section}
\makeatletter
\renewcommand*\l@section{\@dottedtocline{1}{0em}{2.0em}}
\makeatother
\allowdisplaybreaks[1]
\title{The Marvelous Meixner Polynomials:\\[0.2em]
Second Kind, First Class}
\author{Olga V. Holtz}
\date{}
\begin{document}
\maketitle
\begin{abstract}
Meixner's classification asks which exponential generating functions
can generate orthogonal polynomials. We follow the continuous Meixner
family from this question to its moments, measures, and continued
fractions. A beta integral and a bilinear identity establish positivity
and orthogonality; Carlitz's Bernoulli and Euler constructions reveal
further connections. Turning the parameter
into a variable leads to rational orthogonality and Wall's fraction.
We derive its Nevanlinna measure and unite the two spectral directions
in a joint law with Cauchy marginals. Even finite pencils produce
couplings with these exact marginals. Laguerre limits explain the gamma
law governing joint extremes. Along the way, Meixner's scientific life,
combinatorics, tangent numbers, and even zeta values illuminate the
same structure.
\end{abstract}


\normalsize

{\small\setlength{\parskip}{0pt}\tableofcontents}
\clearpage


\section{A question worth asking}\label{sec:opening}

What do the derivatives of the tangent function, a family of orthogonal
polynomials, and a continued fraction have to do with one another? A
three-term recurrence gives a surprisingly economical answer. Fix
$\eta>0$ and $\delta\in\R$, and define monic polynomials by
\begin{equation}\label{eq:recurrence}
 M_{n+1}(x;\delta,\eta)
 =[x-(2n+\eta)\delta]M_n(x;\delta,\eta)
 -n(n+\eta-1)(1+\delta^2)M_{n-1}(x;\delta,\eta),
\end{equation}
with $M_{-1}=0$ and $M_0=1$. The first few are
\begin{align*}
 M_0&=1,\qquad M_1=x-\eta\delta,\\
 M_2&=x^2-2(\eta+1)x\delta+\eta(\eta+1)\delta^2-\eta.
\end{align*}
These are the Meixner polynomials of the second kind, more commonly
encountered today in the Meixner--Pollaczek normalization. Already their
name conceals a small history. Their recurrence conceals rather more.

We shall find that their normalized moments have the exponential
generating function
\begin{equation}\label{eq:intro-mgf}
 H_{\delta,\eta}(u)=(\cos u-\delta\sin u)^{-\eta}.
\end{equation}
At $\delta=0$, $\eta=1$, the moments are the secant numbers. At
$\delta=0$, $\eta=2$, they are the derivatives of $\tan u$, shifted
by one place. Thus classical counting sequences arise from a perfectly
ordinary positive measure. We will see the measure, prove its
orthogonality, and understand why its continued fraction has such
simple coefficients.

Then we will turn the same recurrence sideways. The Meixner-II are
polynomials both in $\delta$ and in $x$. What becomes of orthogonality
when the parameter becomes the variable? The coefficient $1+\delta^2$
hints that the answer should involve rational functions. It does:
the second family is the complementary Romanovski--Routh family, and
its continued fraction belongs to Wall's theory. Better still, both
measures come from one explicit positive kernel. With the correct
normalization, they are the two conditional laws of a joint probability
measure with Cauchy marginals.

The purpose of this paper is to make these connections understandable
without requiring the reader to learn several theories before seeing
what they explain. The order is deliberate. We begin with the
classification problem that led Meixner to the polynomials. Generating
functions then supply the polynomials and their moments together.
Fourier inversion produces the measure, and a product of two generating
functions proves orthogonality. Carlitz's work provides a second point
of entry: start with familiar numbers, treat them as a moment
functional, and recover the polynomials. Paths make the coefficient
arrays visible; finite Jacobi matrices make continued-fraction
convergence precise. Only then do we need the second spectral variable.

The second measure will, in turn, raise two questions of its own.
Does the joint law already exist before the continued fractions are
allowed to go to infinity? And what makes two Cauchy variables with
fixed marginals depend on one another? Finite matrices answer the first
question exactly; a gamma law at the edge answers the second. The same
Laguerre polynomials describe both the directions of the zero curves
and the limiting fluctuations. Two appendices pursue the corresponding
operator and zero geometry in greater detail.

Most of the individual ingredients are classical. We distinguish their
sources from the direct derivations that assemble them here; in
particular, we make no novelty claim for the CRR identification or the
general Wall correspondence. The point is to see how much of the
picture follows from a few calculations, and where a change of
viewpoint reveals something hidden.

We use $(a)_n=a(a+1)\cdots(a+n-1)$, with $(a)_0=1$, and write
$\dd x$ for integration with respect to Lebesgue measure. All formal
power-series identities are explicitly distinguished from analytic
identities. A polynomial sequence is called an \emph{orthogonal
polynomial sequence} when its members have successive degrees, are
pairwise orthogonal, and have nonzero squared norms for the indicated
moment functional; a positive
orthogonality measure is stated separately. These distinctions will
matter in the discussion of Carlitz.


\section{Meixner and the lives of his polynomials}\label{sec:history}

An eponym can make a mathematical object seem to have had a single birth.
The Meixner polynomials have a more interesting history.  They arise from
classification, from limits, from moment problems, and from the search for
useful formulas.  These approaches explain different features of the same
objects.  Keeping them in view will help us understand why changing the
variable in a recurrence can reveal something that its original name hides.

\subsection{A life in mathematical physics}

Josef Meixner (1908--1994) was trained in Munich, where he studied
mathematics, physics, and physical chemistry and completed his doctorate
under Arnold Sommerfeld in 1931.  After assistantships in Munich and
Giessen, he habilitated in 1936/37.  His appointment to Aachen in 1942
was interrupted by military service; he began work there after the war,
became a regular chairholder in 1951, and retired in 1974.
Butzer and Koornwinder give a documented account of his career
\cite[\S2]{butzer2019josef}.

That account also records Meixner's membership in the SA from 1934 and
the NSDAP from 1937.  His postwar employment required a denazification
certificate, and Sommerfeld supplied a supporting letter in December
1945
\cite[\S2]{butzer2019josef}.

The mathematical range of his work is best seen through the problems he
chose.  His dissertation concerned the Green function of the
wave-mechanical Kepler problem \cite{Meixner1933Green}.  A Green function
packages solutions of an equation into a kernel: the choice of solutions,
their normalization, and their behavior at the endpoints all matter at
once.  This is a natural meeting place for mathematical physics and the
theory of special functions.

Meixner's later publications continued that conversation.  His book with
Friedrich Wilhelm Sch\"afke, published in 1954, develops Mathieu and
spheroidal functions and their physical and technical applications
\cite{MeixnerSchaefke1954}.  The organization is revealing: foundations,
the two families of functions, then applications.  Explicit formulas and
structural understanding are meant to support each other.

Another substantial part of his work concerned irreversible processes.
The 1959 \emph{Handbuch der Physik} chapter with H.~G.~Reik develops a
macroscopic thermodynamics extending the equilibrium theory to evolving
continuous media \cite{MeixnerReik1959}.  Its scope includes fluid motion,
elasticity, and media in electromagnetic fields.  Here the unifying
question is how different physical processes fit into a common
mathematical description.  With Fritz Bopp, Meixner also edited
Sommerfeld's \emph{Thermodynamik und Statistik}, published in 1952
\cite{Sommerfeld1952}.  His scientific life thus extends well beyond the
eight-page paper that attached his name to two polynomial families.

\subsection{The question behind the classification}

That paper, received in May 1933 and published in 1934, begins with a
question one can ask before knowing any special-function terminology.
Suppose an exponential generating function has the form
$A(t)e^{xU(t)}$, with $A(0)=1$, $U(0)=0$, and $U'(0)=1$.
When can its polynomials also be orthogonal for a positive measure on
the real line? Hermite and Laguerre polynomials provide affirmative examples.
Bernoulli and Euler polynomials have equally attractive generating
functions, but fail the orthogonality requirement.  The resemblance of
their formulas is therefore a reason to ask the question, not an answer
to it \cite[\S1]{meixner1934orthogonale}.

The starting point was the simpler case $U(t)=t$, now called the Appell
case.  Differentiation then gives $p_n'=np_{n-1}$: every polynomial
lowers to the previous one.  Meixner explicitly credits Salomon Bochner
with participation in the result for this special case, where only the
Hermite family survives.  His acknowledgment should not be enlarged
into a joint authorship of the full classification
\cite[p.~7, footnote]{meixner1934orthogonale}.

Allowing $U$ to vary changes the lowering operator, while orthogonality
still demands a three-term recurrence.  Meixner's decisive step is to
make these two requirements speak to each other.  The resulting
restrictions reduce the problem to a quadratic polynomial and its root
configurations.  The answer consists, up to the usual changes of scale
and origin, of Hermite, Laguerre, Charlier, discrete Meixner, and the
continuous family now called Meixner--Pollaczek
\cite[\S\S4--6]{meixner1934orthogonale}.

This is a classification of \emph{infinite} positive orthogonal systems:
Meixner explicitly assumes a measure with infinitely many points of
increase.  Finite Krawtchouk systems require a change of
the question.  They are not a counterexample to the theorem as he
stated it.  Nor are the two families bearing his name both discrete:
their measures live on a lattice and on the real line, respectively.

Meixner also recovered the moment transform by integrating the
generating function and using orthogonality to the constant polynomial.
Inverting the inner series led to Fourier or Laplace inversion and
the measure itself \cite[\S6]{meixner1934orthogonale}. We carry out
this calculation in Sections~\ref{sec:generating}--\ref{sec:measure}.
The classification tells us when the encoding in the generating
function is compatible with positivity.

\subsection{Discovery, rediscovery, and a neighboring family}

The discrete family had appeared in Ladislav Truksa's 1931 work under
the name \emph{generalized Kummer polynomials}
\cite[Part~III, \S6, pp.~198--199]{truksa1931}.  Truksa approached them
through limiting transitions in hypergeometric orthogonal systems.  He
gave their weight, hypergeometric and Rodrigues formulas, recurrence,
and difference equation.  His viewpoint makes the connection with
probability particularly visible: the weight is a generalized Poisson
law, the negative binomial distribution in present terminology.
Meixner's later classification supplies a different explanation for
their existence.  One construction locates them among limiting
families; the other shows that they are forced by the coexistence of a
generating function and orthogonality.

The continuous family acquired another life in F\'elix Pollaczek's 1950
paper, whose title advertises Hermite and Laguerre polynomials as
limiting cases \cite{Pollaczek1950Limits}.  Pollaczek obtained it by
a limiting passage from his broader family of orthogonal polynomials.
The name Meixner--Pollaczek preserves both routes.  For our purposes the
distinction between the continuous and discrete measures matters more
than the shared part of their names: analytic continuation of formulas
does not automatically continue positivity of a measure.

The generating functions themselves also received a new vocabulary.
Sheffer's 1939 paper develops polynomial sets of type zero
\cite{Sheffer1939}; today we speak of Sheffer sequences.  The
generating-function and lowering-operator descriptions needed for
Meixner's classification were already present in the 1934 paper.
The useful modern distinction is between the broad algebraic class of
Sheffer sequences and its much smaller positive orthogonal subclass.
It is the intersection that makes Meixner's theorem so restrictive.

Finally, Bateman's 1933 polynomials belong to a neighboring story
\cite{bateman1933properties}.  Their route through Legendre polynomials
and Fourier transforms leads to the continuous Hahn family.  Koelink
traces this development through Pasternack, Touchard, and Carlitz
\cite[\S1]{Koelink1996}.  This distinction will be essential when we
read Carlitz: the Bernoulli and Euler sides of his paper need not land
in the same family just because their moment formulas are similar.
History here provides a guide for the mathematics:
compare the recurrence and its measure before identifying the
polynomials by name.


\section{Why only five families?}\label{sec:classification}

Meixner's question can be expressed in one line. Which orthogonal
polynomial sequences have an exponential generating function
\begin{equation}\label{eq:sheffer-general}
 \mathcal G(x;t)=f(t)e^{xg(t)}
 =\sum_{n\ge0}p_n(x)\frac{t^n}{n!},
 \qquad f(0)=1,\quad g(0)=0,\quad g'(0)=1?
\end{equation}
The normalization makes the polynomials monic. Such sequences are now
called \emph{Sheffer sequences}. The exponential in
\eqref{eq:sheffer-general} allows differentiation in $x$ to be
replaced by multiplication by $g(t)$. Orthogonality, on the other
hand, asks multiplication by $x$ to have only three diagonals in the
polynomial basis. The classification comes from making these two
requirements agree.

Indeed, differentiation with respect to $t$ gives
\begin{equation}\label{eq:general-multiplication}
 x\mathcal G=A(t)\mathcal G_t+B(t)\mathcal G,
 \qquad
 A(t)=\frac1{g'(t)},\qquad
 B(t)=-\frac{f'(t)}{f(t)g'(t)}.
\end{equation}
Write $A(t)=\sum A_jt^j$, $B(t)=\sum B_jt^j$, with $A_0=1$.
For $j\ge2$, the coefficient of $p_{n-j}$ in $xp_n$ is
\[
 \frac{n!}{(n-j)!}\bigl[(n-j)A_{j+1}+B_j\bigr].
\]
A three-term recurrence requires this to vanish for every $n\ge j$.
Taking $n=j$ and then $n=j+1$ gives $B_j=A_{j+1}=0$.
Thus necessarily
\begin{equation}\label{eq:quadratic-affine}
 A(t)=1+at+bt^2,\qquad B(t)=\alpha+\beta t.
\end{equation}
Conversely, these forms give exactly
\begin{equation}\label{eq:classification-rec}
 xp_n=p_{n+1}+(an+\alpha)p_n
       +n\bigl[b(n-1)+\beta\bigr]p_{n-1}.
\end{equation}
For a positive measure on the real line the coefficients
$a,b,\alpha,\beta$ are real. For an infinite positive orthogonality
measure, the last coefficient
must be positive for every $n\ge1$. This is equivalent to
$\beta>0$ and $b\ge0$. Favard's theorem gives the converse
\cite{chihara1978introduction}.

We have reduced the problem to a quadratic polynomial. Its possibilities,
up to affine changes of the variable and harmless normalizations, are
\begin{center}
\begin{tabular}{lll}
\toprule
Condition & Family & Support of a standard positive measure\\
\midrule
$b=0$, $a=0$ & Hermite & $\R$\\
$b=0$, $a\ne0$ & Charlier & $\mathbb Z_{\ge0}$\\
$b>0$, $a^2-4b=0$ & Laguerre & $[0,\infty)$\\
$b>0$, $a^2-4b>0$ & Meixner of the first kind & $\mathbb Z_{\ge0}$\\
$b>0$, $a^2-4b<0$ & Meixner of the second kind & $\R$\\
\bottomrule
\end{tabular}
\end{center}
This is the substance of Meixner's five-family classification
\cite{meixner1934orthogonale}. It is a classification of infinite
positive systems; finite systems require separate endpoint conditions.

The case of interest to us has
\[
 a=2\delta,\quad b=1+\delta^2,\quad
 \alpha=\eta\delta,\quad\beta=\eta(1+\delta^2).
\]
Its discriminant is always $-4$. More generally, if $a^2<4b$, put
\[
 \omega=\sqrt{b-a^2/4},\qquad
 \delta=\frac{a}{2\omega},\qquad\eta=\frac\beta b.
\]
The monic rescaling
\[
 q_n(X)=\omega^{-n}p_n(\omega X+\alpha-\eta a/2)
\]
reduces \eqref{eq:classification-rec} to \eqref{eq:recurrence} in $X$.
The quadratic factor $1+\delta^2$ is therefore a consequence of a
particularly convenient normalization of the negative-discriminant
case. It is also the factor that will later turn the parameter
recurrence into an $R_{\mathrm{II}}$ recurrence.

The calculation above is worth remembering. A rather grand
classification theorem begins with the observation that a linear
function of $n$ cannot vanish at every integer unless both of its
coefficients vanish. The special functions enter only after this
restriction has done its work.


\section{One recurrence, two generating functions}\label{sec:generating}

\subsection{The polynomials}

Let
\[
 G(x;t)=\sum_{n\ge0}M_n(x;\delta,\eta)\frac{t^n}{n!},
 \qquad D(t)=1+2\delta t+(1+\delta^2)t^2.
\]
Multiply \eqref{eq:recurrence} by $t^n/n!$ and sum. The only term
whose index deserves a second look is
\[
 \sum_{n\ge1}n(n+\eta-1)M_{n-1}\frac{t^n}{n!}
 =t\sum_{m\ge0}(m+\eta)M_m\frac{t^m}{m!}
 =t^2G_t+\eta tG.
\]
Consequently,
\begin{equation}\label{eq:G-ode}
 D(t)G_t=[x-\eta\delta-\eta(1+\delta^2)t]G,
 \qquad G(x;0)=1.
\end{equation}
This is an ordinary first-order differential equation in $t$; $x$ is
a parameter. Its logarithmic derivative separates immediately:
\[
 \frac{G_t}{G}=\frac{x}{D(t)}-\frac\eta2\frac{D'(t)}{D(t)}.
\]
We obtain
\begin{equation}\label{eq:G}
 G(x;t)=D(t)^{-\eta/2}e^{xg(t)},\qquad
 g(t)=\int_0^t\frac{dv}{D(v)}
     =\arctan\frac{t}{1+\delta t}.
\end{equation}
The last expression means the branch analytic at $0$, with value $0$
there. Equivalently, the two scalar equations in the Sheffer form
$G=f e^{xg}$ are
\[
 g'=D^{-1},\qquad f'/f=-\eta D'/(2D).
\]
Thus the generating function is not a guess. The recurrence
has already specified its two factors.

\subsection{Inversion, and what the word means}

Because the $M_n$ are monic, there are unique numbers $q_{n,k}$ such
that
\begin{equation}\label{eq:connection}
 x^n=\sum_{k=0}^nq_{n,k}M_k(x;\delta,\eta).
\end{equation}
The \emph{inverse polynomials} are
$Q_n(y)=\sum_{k=0}^nq_{n,k}y^k$. This is inversion of a triangular
coefficient matrix, not inversion of a polynomial as a function.

The compositional inverse of $g$ in \eqref{eq:G} is
\begin{equation}\label{eq:tau}
 \tau(u)=g^{-1}(u)=\frac{\tan u}{1-\delta\tan u}
                      =\frac{\sin u}{\cos u-\delta\sin u}.
\end{equation}
Substituting $t=\tau(u)$ in $G=f e^{xg}$ gives
\[
 e^{xu}=\frac1{f(\tau(u))}
       \sum_{k\ge0}M_k(x)\frac{\tau(u)^k}{k!}.
\]
Comparison with \eqref{eq:connection} proves the Sheffer inversion
formula in this case:
\begin{equation}\label{eq:inverse-egf}
 \sum_{n\ge0}Q_n(y)\frac{u^n}{n!}
 =H(u)e^{y\tau(u)},\qquad
 H(u)=\frac1{f(\tau(u))}
     =(\cos u-\delta\sin u)^{-\eta}.
\end{equation}
We used the simple identity $D(\tau(u))=(\cos u-\delta\sin u)^{-2}$.
In particular, the exponent in the prefactor of the inverse generating
function is $\eta$, not $\eta/2$.

This is one useful meaning of \emph{umbral calculus}: a linear map
sends the monomial basis to a polynomial basis, and substitution and
composition organize the map and its inverse. There is also a lowering
operator. If $D_x=d/dx$, then, as a formal differential operator on
polynomials,
\begin{equation}\label{eq:lowering-sheffer}
 \tau(D_x)M_n=nM_{n-1}.
\end{equation}
Indeed, $\tau(D_x)G=\tau(g(t))G=tG$. No convergence question is
hidden here: on a polynomial only finitely many derivatives survive.
This operator language belongs to the same tradition as Sheffer's
classification by lowering operators \cite{Sheffer1939}.

The inverse sequence is again Sheffer, but need not be an ordinary
orthogonal sequence. For example, when $\delta=0$, $\eta=1$,
\[
 Q_0=1,\quad Q_1=y,\quad Q_2=y^2+1,\quad
 Q_3=y^3+5y,\quad Q_4=y^4+14y^2+5.
\]
Here $yQ_3=Q_4-9Q_2+4Q_0$, which rules out a three-term recurrence.
Triangular inversion preserves the Sheffer structure, not automatically
orthogonality.

\subsection{The moments are already in the inverse matrix}

Define a linear functional $\mathcal L$ on polynomials by
$\mathcal L(M_0)=1$ and $\mathcal L(M_n)=0$ for $n\ge1$. At
this point the definition is algebraic; the measure will come next.
Applying $\mathcal L$ to \eqref{eq:connection} gives
\[
 \mu_n:=\mathcal L(x^n)=q_{n,0}=Q_n(0).
\]
Putting $y=0$ in \eqref{eq:inverse-egf}, we get all moments at once:
\begin{equation}\label{eq:H}
 \sum_{n\ge0}\mu_n(\delta,\eta)\frac{u^n}{n!}
 =H_{\delta,\eta}(u)=(\cos u-\delta\sin u)^{-\eta}.
\end{equation}
Equivalently, $\mathcal L(G(x;t))=1$ says
$f(t)H(g(t))=1$. The inversion formula and the moment calculation are
two readings of the same identity.

For later checks, the first moments are
\begin{align*}
 \mu_0&=1,\qquad\mu_1=\eta\delta,\\
 \mu_2&=\eta+\eta(\eta+1)\delta^2,\\
 \mu_3&=\eta(3\eta+2)\delta
             +\eta(\eta+1)(\eta+2)\delta^3.
\end{align*}
A convenient recurrence follows by differentiating \eqref{eq:H}:
\begin{equation}\label{eq:moment-derivative}
 \mu_{n+1}=(1+\delta^2)\partial_\delta\mu_n+\eta\delta\mu_n,
 \qquad\mu_0=1.
\end{equation}
For $\delta=0$, the generating function is $\sec^\eta u$. It is
important not to confuse this with $G(0;t)=(1+t^2)^{-\eta/2}$.
The former generates moments; the latter generates the values of the
orthogonal polynomials at the origin. They answer different questions.


\section{A beta integral finds the measure}\label{sec:measure}

The moments have appeared algebraically. We now construct their positive
measure and prove orthogonality directly. The main calculation is a
Fourier transform that becomes a beta integral after one substitution.

\subsection{The symmetric density and its exponential tilt}

For $\eta>0$, set
\begin{equation}\label{eq:Ceta}
 C_\eta(x)=\frac{2^{\eta-2}}{\pi\Gamma(\eta)}
       \left|\Gamma\!\left(\frac{\eta+ix}{2}\right)\right|^2.
\end{equation}
To discover this expression, put $v=e^{2u}$ in the Fourier integral:
\begin{align}\label{eq:beta-Fourier}
 \int_{\R}e^{-ixu}\operatorname{sech}^{\eta}u\,du
 &=2^{\eta-1}\int_0^\infty
       \frac{v^{(\eta-ix)/2-1}}{(1+v)^\eta}\,dv\notag\\
 &=\frac{2^{\eta-1}}{\Gamma(\eta)}
       \Gamma\!\left(\frac{\eta-ix}{2}\right)
       \Gamma\!\left(\frac{\eta+ix}{2}\right).
\end{align}
Fourier inversion therefore gives
\begin{equation}\label{eq:symmetric-characteristic}
 \int_{\R}e^{iux}C_\eta(x)\,dx
 =\operatorname{sech}^{\eta}u.
\end{equation}
All integrals are ordinary absolutely convergent integrals. For example,
Stirling's formula gives
\begin{equation}\label{eq:Ceta-tail}
 C_\eta(x)\sim \frac{|x|^{\eta-1}}{\Gamma(\eta)}
                    e^{-\pi|x|/2}\qquad(|x|\to\infty).
\end{equation}
The density is positive, and its integral is one by
\eqref{eq:symmetric-characteristic} at $u=0$.

Write $\theta=\arctan\delta$, so that
$-\pi/2<\theta<\pi/2$. Exponentially tilting this symmetric probability
gives our central kernel:
\begin{equation}\label{eq:kernel}
 \boxed{\displaystyle
 K_\eta(x,\delta)
 =C_\eta(x)(1+\delta^2)^{-\eta/2}e^{x\arctan\delta}.}
\end{equation}
Indeed, the exponential tails allow analytic continuation of
\eqref{eq:symmetric-characteristic} to $|\Im u|<\pi/2$, and hence
$\int e^{\theta x}C_\eta(x)\,dx=\cos^{-\eta}\theta$.
Thus
\begin{equation}\label{eq:mu-measure}
 d\mu_\delta(x)=K_\eta(x,\delta)\,dx
\end{equation}
is a probability measure. Its moment generating function is
\begin{align}\label{eq:measure-mgf}
 \int_{\R}e^{ux}\,d\mu_\delta(x)
 &=\frac{\cos^\eta\theta}{\cos^\eta(\theta+u)}\notag\\
 &=(\cos u-\delta\sin u)^{-\eta},
 \qquad -\frac\pi2<\theta+u<\frac\pi2.
\end{align}
This agrees with the formal series \eqref{eq:H}. Our algebraic moment
functional has found its positive realization.

The measure is also unique among positive measures with these moments.
Here is a short reason. For sufficiently small $r>0$,
$\sum_{n\ge0}\mu_{2n}r^{2n}/(2n)!<\infty$.
If another positive measure has the same moments, monotone convergence
applied to the series of $\cosh(rx)$ gives the same finite integral.
Its moment generating function is consequently analytic near zero,
with the same Taylor coefficients as \eqref{eq:measure-mgf}.
Equality extends through their common vertical strip, in particular
along the imaginary axis. Uniqueness of the Fourier transform then
gives the same measure.
This observation will also settle the convergence of our Jacobi fraction.

\subsection{Orthogonality in one line, once the line is understood}

Let $G(x;t)=D(t)^{-\eta/2}e^{xg(t)}$ be the polynomial generating
function from \eqref{eq:G}, with
$D(t)=1+2\delta t+(1+\delta^2)t^2$.
For small $s,t$, equation \eqref{eq:measure-mgf} gives
\[
 \int G(x;s)G(x;t)\,d\mu_\delta(x)
 =D(s)^{-\eta/2}D(t)^{-\eta/2}
   [\cos(g(s)+g(t))-\delta\sin(g(s)+g(t))]^{-\eta}.
\]
The trigonometry has been arranged for us by the recurrence:
\[
 \cos g(t)=\frac{1+\delta t}{\sqrt{D(t)}},\qquad
 \sin g(t)=\frac{t}{\sqrt{D(t)}}.
\]
The addition formulas therefore reduce the bracket to
\[
 \frac{1-(1+\delta^2)st}{\sqrt{D(s)D(t)}}.
\]
Everything cancels except
\begin{equation}\label{eq:bilinear-G}
 \int_{\R}G(x;s)G(x;t)\,d\mu_\delta(x)
   =\bigl(1-(1+\delta^2)st\bigr)^{-\eta}.
\end{equation}
Comparing coefficients proves both orthogonality and the norms:
\begin{equation}\label{eq:orthogonality}
 \boxed{\displaystyle
 \int_{\R}M_m(x;\delta,\eta)M_n(x;\delta,\eta)
           K_\eta(x,\delta)\,dx
 =n!(\eta)_n(1+\delta^2)^n\,\delta_{mn}.}
\end{equation}
The coefficient comparison is justified analytically in a neighborhood
of $(s,t)=(0,0)$ by the exponential tails; it may equally be read
coefficient by coefficient after the moments have been identified.
No measure has been guessed from the name of the polynomials.

For comparison with the conventional notation, put
\begin{equation}\label{eq:MP-dictionary}
 \lambda=\frac\eta2,\qquad
 \phi=\frac\pi2+\arctan\delta,\qquad \delta=-\cot\phi.
\end{equation}
Then the standard Meixner--Pollaczek polynomials satisfy
\begin{equation}\label{eq:MP-normalization}
 M_n(x;\delta,\eta)
 =\frac{n!}{\sin^n\phi}
       P_n^{(\lambda)}(x/2;\phi).
\end{equation}
The identity follows either from their generating function or by
matching the recurrence and the initial values
\cite[\S9.7]{KoekoekLeskySwarttouw2010}. The angle convention fixes
the sign of the exponential tilt in \eqref{eq:kernel}.

\subsection{The two parameters have simple jobs}

At fixed $\delta$, multiplying the moment generating functions gives
\begin{equation}\label{eq:convolution-general}
 \mu_{\delta,\eta_1}*\mu_{\delta,\eta_2}
       =\mu_{\delta,\eta_1+\eta_2},\qquad \eta_1,\eta_2>0.
\end{equation}
Thus $\eta$ adds independent summands. The parameter $\delta$, or
equivalently $\theta$, tilts the symmetric law toward one end of the
line. In particular,
\[
 \mathbb E[X]=\eta\delta,\qquad
 \operatorname{Var}(X)=\eta(1+\delta^2).
\]
There is a parallel polynomial addition formula. Multiplying
$G_{\eta_1}(x;t)$ and $G_{\eta_2}(y;t)$ and comparing coefficients gives
\begin{equation}\label{eq:polynomial-addition}
 M_n(x+y;\delta,\eta_1+\eta_2)
 =\sum_{k=0}^n\binom nk
   M_k(x;\delta,\eta_1)M_{n-k}(y;\delta,\eta_2).
\end{equation}
The probability law and the polynomial sequence express the same
additive structure. Carlitz's higher-order Euler numbers will make
this correspondence especially concrete.

\subsection{The shape parameter becomes time}\label{sec:meixner-time}

There is another way to read the addition law: let the shape parameter
measure elapsed time. Fix $\delta$, and write
\[
 \mu_t^{(\delta)}(dx)=K_t(x,\delta)\,dx,\qquad t>0.
\]
These probabilities form a convolution semigroup by
\eqref{eq:convolution-general}. Their characteristic functions tend to
one as $t\downarrow0$, so the limiting law is a unit mass at zero.
They therefore describe a process $(X_t)_{t\ge0}$ starting at zero,
with stationary independent increments and
$X_t-X_s\sim\mu_{t-s}^{(\delta)}$ for $s<t$. This is the
\emph{Meixner process}. For our calculations, only its joint laws at
finitely many times are needed: prescribe independent increments with
these densities; the convolution identity makes the prescriptions
consistent. The connection with orthogonal polynomial martingales is
classical \cite{SchoutensTeugels1998}, but its proof here is already
contained in the addition formula.

Abbreviate $M_n^{(t)}(x)=M_n(x;\delta,t)$ and set
\begin{equation}\label{eq:time-norms}
 h_n(t)=n!(t)_n(1+\delta^2)^n.
\end{equation}

\begin{proposition}\label{prop:research-two-martingales}
For $0<s<t$ and $n\ge0$,
\begin{align}
 \E[M_n^{(t)}(X_t)\mid\mathcal F_s]
   &=M_n^{(s)}(X_s),\label{eq:research-forward-martingale}\\
 \E[M_n^{(s)}(X_s)\mid X_t]
   &=\frac{(s)_n}{(t)_n}M_n^{(t)}(X_t),
     \label{eq:research-reverse-martingale}
\end{align}
where $\mathcal F_s$ records the process up to time $s$.
\end{proposition}

\begin{proof}
Write $X_t=X_s+Y$, with $Y$ independent of $\mathcal F_s$ and
having law $\mu_{t-s}^{(\delta)}$. Expand the polynomial by
\eqref{eq:polynomial-addition}. Orthogonality to the constant says
$\E[M_j(Y;\delta,t-s)]=0$ for $j>0$, so conditional expectation
leaves exactly the term $M_n^{(s)}(X_s)$.

For the other direction, the identity just proved and orthogonality give
\[
 \E[M_n^{(s)}(X_s)M_m^{(t)}(X_t)]
 =\E[M_n^{(s)}(X_s)M_m^{(s)}(X_s)]
 =\delta_{nm}h_n(s).
\]
Thus the orthogonal expansion of
$\E[M_n^{(s)}(X_s)\mid X_t]$ has only its $n$th term, with
coefficient $h_n(s)/h_n(t)=(s)_n/(t)_n$.

To justify that expansion, the polynomials are complete in
$L^2(\mu_t^{(\delta)})$. Indeed, if $f$ is orthogonal to every
polynomial, Cauchy--Schwarz and the exponential tails make
$\int f(x)e^{wx}\,\mu_t^{(\delta)}(dx)$ analytic in a strip about
the imaginary axis. All its derivatives at zero vanish. The identity
theorem, followed by uniqueness of the Fourier transform, gives
$f=0$ almost everywhere. Conditional expectation belongs to this
$L^2$ space by Jensen's inequality, so its expansion identifies it
uniquely.
\end{proof}

The first identity says that the polynomial's present value is its
conditional future average: it is a martingale. The second gives the
reverse statement after division by $(t)_n$. More precisely, if
$\mathcal F^t=\sigma(X_r:r\ge t)$ records the future, independence of
the increments after $t$ shows that
$M_n^{(t)}(X_t)/(t)_n$ is a reverse martingale with respect to these
decreasing sigma-fields.

The first two degrees already give useful conditional moments. Since
\[
 M_1^{(t)}(x)=x-t\delta,\qquad
 M_2^{(t)}(x)=x^2-2(t+1)\delta x+t(t+1)\delta^2-t,
\]
substitution in \eqref{eq:research-reverse-martingale} yields
\begin{equation}\label{eq:research-bridge-moments}
 \boxed{\displaystyle
 \E[X_s\mid X_t=x]=\frac{s}{t}x,\qquad
 \operatorname{Var}(X_s\mid X_t=x)
 =\frac{s(t-s)}{t+1}\left(1+\frac{x^2}{t^2}\right).}
\end{equation}
The tilt has disappeared. In the next section we will see why the
entire conditional density shares this property, and why Carlitz's
Bernoulli construction appears in it. These conditional moments are
also recorded in \cite{DevroyeHill2025}. Appendix~\ref{app:conditioning}
uses the same two polynomial identities to determine exactly how much
information survives conditioning.


\section{Carlitz's route: numbers first, polynomials later}\label{sec:carlitz}

Carlitz enters this story from the opposite direction. He starts with familiar
numbers---Bernoulli numbers, Euler numbers, and their relatives---and looks for
polynomials that make those numbers into moments. His 1959 paper
\cite{carlitz1959bernoulli} moves repeatedly between a symbolic calculation, a
three-term recurrence, and an integral or sum. Reading it through Meixner's
classification reveals how much of the same landscape these different starting
points uncover.

This is more than a change of notation. A moment sequence can look rather
unpromising when inspected term by term, yet become transparent in a basis
adapted to it. Carlitz repeatedly finds that basis by replacing powers with
binomial coefficients. The resulting cancellation is elementary; its
consequences include orthogonality, explicit norms, and eventually a measure.
We follow that thread through the parts of his paper that meet the Meixner
family, and then briefly examine the Bernoulli branch that goes beyond it.

The symbolic language deserves a brief explanation. If numbers $a_n$ are given,
the notation $a^n=a_n$ means that the linear functional
$\mathcal A[u^n]=a_n$ is to be applied after expanding a polynomial. It does not
assert that the $a_n$ are powers of one number. Likewise,
$e^{at}=\sum a_nt^n/n!$ is shorthand for $\mathcal A[e^{ut}]$.
Such a functional need not be positive. That distinction will explain both the
imaginary changes of variable below and the signs in Carlitz's formulas.

\subsection{Eulerian moments: one calculation, two kinds of Meixner polynomials}

Carlitz's ``Eulerian'' numbers are defined by
\begin{equation}\label{eq:car-eulerian-functional}
 \mathcal R_\lambda[e^{zt}]
 =\sum_{j\ge0}R_j(\lambda)\frac{t^j}{j!}
 =\frac{1-\lambda}{e^t-\lambda},\qquad \lambda\ne1.
\end{equation}
The associated polynomials have generating function
\begin{equation}\label{eq:car-A-generating}
 \sum_{n\ge0}A_n^{(\lambda)}(z)t^n
 =\frac{(1-\lambda t)^z}{(1-t)^{z+1}}.
\end{equation}
Here the series is ordinary in $t$; the leading coefficient of
$A_n^{(\lambda)}$ is $(1-\lambda)^n/n!$. These are
\cite[(1.17), (1.19)]{carlitz1959bernoulli}; the parameter $\lambda$
here is Carlitz's parameter, not our shape parameter $\eta$.

Carlitz's cancellation is particularly pleasant. Substituting
$t=\log(1+v)$ in \eqref{eq:car-eulerian-functional} gives
\[
 \mathcal R_\lambda\!\left[\binom zr\right]
 =(\lambda-1)^{-r}.
\]
On the other hand, expanding \eqref{eq:car-A-generating} gives
\begin{equation}\label{eq:car-A-binomial}
 A_n^{(\lambda)}(z)
 =\sum_{r=0}^n(1-\lambda)^r\binom nr\binom zr
 ={}_2F_1(-n,-z;1;1-\lambda).
\end{equation}
Applying $\mathcal R_\lambda$ leaves
$\sum_r(-1)^r\binom nr=0$ for $n>0$.
This is the essential step in \cite[(5.1)--(5.4)]{carlitz1959bernoulli}.
Differentiating the generating function supplies the recurrence
\begin{equation}\label{eq:car-A-recurrence}
 (n+1)A_{n+1}^{(\lambda)}(z)
 =\bigl[(1-\lambda)z+1+n(1+\lambda)\bigr]A_n^{(\lambda)}(z)
   -n\lambda A_{n-1}^{(\lambda)}(z).
\end{equation}
Multiplication by $z$ can lower the polynomial index by at most one. Therefore
$\mathcal R_\lambda[z^rA_n^{(\lambda)}]=0$ for $r<n$: fewer than $n$
multiplications cannot reach the constant polynomial. At the first possible
step there is just one contribution: every multiplication must lower the
index. Solving \eqref{eq:car-A-recurrence} for $zA_n^{(\lambda)}$ shows
that the lowering coefficient is $n\lambda/(1-\lambda)$, so
\[
 \mathcal R_\lambda[z^n A_n^{(\lambda)}]
 =\prod_{j=1}^n\frac{j\lambda}{1-\lambda}
 =\frac{n!\lambda^n}{(1-\lambda)^n}.
\]
Multiplying by the leading coefficient of $A_n^{(\lambda)}$ gives
\begin{equation}\label{eq:car-A-orthogonality}
 \mathcal R_\lambda[A_m^{(\lambda)}A_n^{(\lambda)}]
 =\lambda^n\delta_{mn}.
\end{equation}
This illustrates a useful feature of Carlitz's method. Once the average of each
nonconstant polynomial is zero, the recurrence does the remaining work.

For $\lambda>1$, the functional is already a positive probability law:
\begin{equation}\label{eq:car-geometric}
 \mathcal R_\lambda[p]
 =(\lambda-1)\sum_{j=0}^\infty p(j)\lambda^{-j-1}.
\end{equation}
Indeed, summing the geometric series for $p(z)=e^{zt}$ gives
\eqref{eq:car-eulerian-functional}. Equation
\eqref{eq:car-A-binomial} now identifies $A_n^{(\lambda)}$ with the
Meixner polynomials of the first kind having shape parameter $1$ and
$c=1/\lambda$. Carlitz gives this discrete orthogonality in
\cite[\S8]{carlitz1959bernoulli}.

The same functional also has the contour representation
\begin{equation}\label{eq:car-contour}
 \mathcal R_\lambda[p]
 =\frac{1-\lambda}{2i\lambda}
  \int_{-1/2-i\infty}^{-1/2+i\infty}
  \frac{p(z)(-\lambda)^{-z}}{\sin\pi z}\,dz,
\end{equation}
for $\lambda\in\mathbb C\setminus[0,\infty)$, with the principal logarithm of
$-\lambda$. This is an absolutely convergent integral for every polynomial:
on a vertical line the reciprocal sine decays like $e^{-\pi|\operatorname{Im}z|}$,
whereas the exponential factor contributes
$e^{\arg(-\lambda)\operatorname{Im}z}$, with $|\arg(-\lambda)|<\pi$.
Carlitz proves the identity by shifting the contour by one unit. The residue
fixes the mass, and the shift gives
\[
 \mathcal R_\lambda[p(z+1)-\lambda p(z)]
 =(1-\lambda)p(0),
\]
which also follows directly from \eqref{eq:car-eulerian-functional}.
For $p(z)=z^n$, the coefficient of the unknown $n$th moment is
$1-\lambda$, and all remaining moments have smaller index. The identity
therefore determines the functional inductively, including its mass when
$n=0$. This is the moment argument behind
\cite[(7.1)--(7.5)]{carlitz1959bernoulli}.

Now put
\begin{equation}\label{eq:car-parameter-change}
 \lambda=\frac{\delta-i}{\delta+i},\qquad
 z=\frac{-1+ix}{2},\qquad \delta\in\mathbb R.
\end{equation}
Since $1-\lambda=2i/(\delta+i)$ and $1+\lambda=2\delta/(\delta+i)$,
the recurrence \eqref{eq:car-A-recurrence} changes into our monic
Meixner recurrence with $\eta=1$. The initial polynomial becomes
$M_1=x-\delta$, fixing the normalization, and we obtain
\begin{equation}\label{eq:car-to-Meixner}
 \boxed{\displaystyle
 M_n(x;\delta,1)
 =n!(-\delta-i)^n
 A_n^{((\delta-i)/(\delta+i))}\!\left(\frac{-1+ix}{2}\right).}
\end{equation}
The contour becomes the real $x$-axis. If $\theta=\arctan\delta$, then
$-\lambda=e^{2i\theta}$, and its weight becomes
\[
 \frac{\cos\theta\,e^{\theta x}}{2\cosh(\pi x/2)}\,dx
 =K_1(x,\delta)\,dx.
\]
Thus Carlitz's contour construction is exactly the positive Meixner-II measure
at $\eta=1$. Even the apparently complex norm in
\eqref{eq:car-A-orthogonality} becomes our positive norm, because
$(-\delta-i)^2\lambda=1+\delta^2$.
The discrete and continuous Meixner families have appeared from the same
symbolic calculation, with different realizations of its functional.
In particular, Carlitz's connection with the continuous family is not
confined to the symmetric case $\delta=0$. His Eulerian parameter already
contains every real tilt $\delta$ at shape $1$. The complex quantities in
\eqref{eq:car-to-Meixner} are coordinates for real polynomials and a positive
measure; the recurrence and the transformed weight make that reality visible.

\subsection{Higher Euler numbers and the meaning of the imaginary rotation}

Carlitz next considers the signed Euler numbers of order $\eta$, which we write
as the moments of $\mathcal E_\eta$:
\begin{equation}\label{eq:car-Euler-umbra}
 \mathcal E_\eta[e^{ut}]
 =\sum_{j\ge0}E_j^{(\eta)}\frac{t^j}{j!}
 =\operatorname{sech}^{\eta}t.
\end{equation}
In particular, $E_2^{(\eta)}=-\eta$. For $\eta>0$ this cannot be the
second moment of a positive measure on the real line. The umbral identity is
nevertheless perfectly meaningful, and positivity will emerge after rotation.

His monic polynomials $f_n^{(\eta)}$ satisfy
\begin{equation}\label{eq:car-f-recurrence}
 f_{n+1}^{(\eta)}(u)
 =u f_n^{(\eta)}(u)+n(n+\eta-1)f_{n-1}^{(\eta)}(u),
 \qquad f_0^{(\eta)}=1,\quad f_1^{(\eta)}=u.
\end{equation}
The plus sign is intentional. If
$\Phi_\eta(u;t)=\sum f_n^{(\eta)}(u)t^n/n!$, the recurrence gives the
first-order equation
\[
 (1-t^2)\partial_t\Phi_\eta=(u+\eta t)\Phi_\eta.
\]
Integration with $\Phi_\eta(u;0)=1$ yields
\begin{equation}\label{eq:car-f-generating}
 \Phi_\eta(u;t)
 =(1-t^2)^{-\eta/2}e^{u\operatorname{arctanh}t}
 =\frac{(1+t)^{(u-\eta)/2}}{(1-t)^{(u+\eta)/2}}.
\end{equation}
This is Carlitz's derivation in \cite[\S9, (9.3)--(9.7)]{carlitz1959bernoulli}.

In his proof of orthogonality, another binomial identity first makes the
nonconstant averages vanish, and the recurrence again propagates the result.
The mechanism can be compressed into one generating-function calculation:
\begin{align}\label{eq:car-bilinear}
 \mathcal E_\eta[\Phi_\eta(u;s)\Phi_\eta(u;t)]
 &=(1-s^2)^{-\eta/2}(1-t^2)^{-\eta/2}
   \operatorname{sech}^{\eta}
       (\operatorname{arctanh}s+\operatorname{arctanh}t)\notag\\
 &=(1+st)^{-\eta}.
\end{align}
Coefficient comparison gives the exact signed norm
\begin{equation}\label{eq:car-signed-norm}
 \mathcal E_\eta[f_m^{(\eta)}f_n^{(\eta)}]
 =(-1)^n n!(\eta)_n\delta_{mn}.
\end{equation}
This argument takes place in formal power series and uses no positivity.
There is no convergence question hidden in the coefficient comparison:
each coefficient uses only finitely many moments. The hyperbolic addition
formula is what removes all mixed terms except those with equal polynomial
indices. It is the same mechanism as our earlier bilinear calculation for
$M_n$, now expressed in the signed Euler convention.

The conversion to our polynomials is now immediate:
\begin{equation}\label{eq:car-imaginary-rotation}
 M_n(x;0,\eta)=i^{-n}f_n^{(\eta)}(ix),\qquad
 \mu_n(0,\eta)=i^{-n}E_n^{(\eta)}.
\end{equation}
The imaginary rotation changes the sign of the recurrence coefficient and
turns the signed norms into $n!(\eta)_n$. In probabilistic terms,
$\mathcal E_\eta[p(u)]=\int p(ix)K_\eta(x,0)\,dx$.
This explains the symbolic substitution without pretending that the signed
Euler numbers themselves are moments of a positive real measure.
\footnote{Our formulas use monic polynomials and the variance-one density
$p_1(x)=1/(2\cosh(\pi x/2))$. Thus the symbolic norm has the factor
$(-1)^n$, and the positive norm at degree one is $\eta$.
These elementary checks fix the sign and scale factors in comparing
\cite[(7.6), (9.15)--(9.16)]{carlitz1959bernoulli}; the printed formulas
there require the corresponding corrections.}

Equivalently, the ordinary moment generating function of the positive law is
$\sec^{\eta}t$ near $t=0$, while its characteristic function is
$\operatorname{sech}^{\eta}t$. This distinction matters: Carlitz's generating
function becomes a characteristic function on the real line, not the
generating function of its ordinary real moments. The factors $i^n$ in
\eqref{eq:car-imaginary-rotation} say precisely this.

For integer $\eta=k\ge1$, Carlitz realizes the measure as a convolution of
$k$ hyperbolic-secant laws. In our scale,
\begin{equation}\label{eq:car-convolution}
 p_1(x)=\frac1{2\cosh(\pi x/2)},\qquad
 p_k=p_1^{*k},\qquad
 \int_{\mathbb R}e^{itx}p_k(x)\,dx
 =\operatorname{sech}^k t.
\end{equation}
Our gamma density identifies $p_k(x)=K_k(x,0)$ and extends this positive
realization to every real $\eta>0$, as in
\eqref{eq:convolution-general}. Carlitz's constructions illuminate both
parameters: the order of the Euler numbers controls convolution, while
the continuous Eulerian parameter realizes the exponential tilt in
\eqref{eq:kernel}.

\subsection{The Bernoulli branch and its intersection with Meixner}

The first part of Carlitz's paper extends Touchard's work on Bernoulli moments.
For a new parameter $a$, define
\begin{equation}\label{eq:car-Bernoulli-beta}
 \beta_j(a)=\frac{B_{j+1}(a)-B_{j+1}}{(j+1)a},\qquad
 \sum_{j\ge0}\beta_j(a)\frac{t^j}{j!}
 =\frac{e^{at}-1}{a(e^t-1)}.
\end{equation}
The limit $a\to0$ gives the Bernoulli numbers. At a positive integer $a=N$,
Faulhaber's formula says that $\beta_j(N)$ is the average of
$0^j,1^j,\ldots,(N-1)^j$. Thus the same expression interpolates a familiar
finite probability law and the Bernoulli moment functional.

Carlitz identifies the corresponding polynomials with the
Bateman--Pasternack family. In his normalization their recurrence is
\begin{equation}\label{eq:car-Omega}
 \Omega_{n+1}^{(a)}(z)
 =(2z+1-a)\Omega_n^{(a)}(z)
 +\frac{n^2(n^2-a^2)}{4n^2-1}\Omega_{n-1}^{(a)}(z).
\end{equation}
The finite support leaves a visible trace: at $a=N$ the recurrence coefficient
vanishes when $n=N$. The nondegenerate orthogonal family then has exactly
$N$ members, corresponding to the $N$ support points of the uniform law.
This is already a different phenomenon from the infinite positive Meixner
families considered above.
After the change $z=(a-1+ix)/2$, the monic polynomials
$Q_n^{(a)}(x)=i^{-n}\Omega_n^{(a)}((a-1+ix)/2)$ satisfy
\[
 Q_{n+1}^{(a)}(x)
 =xQ_n^{(a)}(x)
  -\frac{n^2(n^2-a^2)}{4n^2-1}Q_{n-1}^{(a)}(x).
\]
For $0<a<1$, Carlitz's contour integral
\cite[(1.21), (6.1)]{carlitz1959bernoulli} becomes the positive density
\begin{equation}\label{eq:car-Bernoulli-density}
 \frac{\sin\pi a}{2a}\,
 \frac{dx}{\cos\pi a+\cosh\pi x}.
\end{equation}
This is a continuous Hahn family, neighboring Meixner in the larger
hypergeometric picture \cite{Koelink1996}. Its recurrence coefficient is generally different
from the quadratic Meixner coefficient. At $a=1/2$, however, it reduces to
$n^2/4$, and
\[
 Q_n^{(1/2)}(x)=2^{-n}M_n(2x;0,1).
\]
The density then becomes $1/\cosh\pi x$. This is the intersection between
the Bernoulli and Euler constructions that Carlitz develops in
\cite[\S\S3--4]{carlitz1959bernoulli}. It explains why both sorts of numbers
occur in his paper without requiring that all of its polynomial families be
Meixner families.

\subsection{The Bernoulli density is a Meixner connection}\label{sec:research-bridges}

The neighboring continuous Hahn family has an even closer connection to
Meixner: it appears when we condition two independent Meixner summands
on their sum. Let $(X_t)$ be the process of
Section~\ref{sec:meixner-time}, with fixed tilt $\delta$. For $0<s<t$,
independence of $X_s$ and $X_t-X_s$ gives the conditional density
\begin{equation}\label{eq:research-bridge-density}
 b_{s,t;x}(y)
 =\frac{K_s(y,\delta)K_{t-s}(x-y,\delta)}{K_t(x,\delta)}
 =\boxed{\frac{C_s(y)C_{t-s}(x-y)}{C_t(x)}}.
\end{equation}
The exponential tilts cancel, as do their normalizing factors. Thus
the whole conditional law of $X_s$ given $X_t=x$ is independent of
$\delta$, explaining the cancellation in
\eqref{eq:research-bridge-moments}. The convolution identity gives its
normalization. Since the denominator is positive, this formula supplies
a continuous version of the conditional density for every real endpoint
$x$, including the endpoint zero.

With $r=t-s$ and $z=y/2$, the density in $z$ is proportional to
\[
 \left|\Gamma\!\left(\frac{s}{2}+iz\right)
        \Gamma\!\left(\frac{r-ix}{2}+iz\right)\right|^2.
\]
This is the continuous Hahn weight with parameters
$(s/2,(r-ix)/2,s/2,(r+ix)/2)$. The same disintegration occurs in
Koelink and Van der Jeugt's positive discrete tensor-product construction
\cite[Proposition~3.3 and its proof]{koelink1998convolutions}:
continuous Hahn polynomials resolve the two summands at fixed sum.
In a statistical setting, Devroye and Hill call the conditional law the
\emph{betaized Meixner--Morris distribution}
\cite[\S11]{DevroyeHill2025}. Formula
\eqref{eq:research-bridge-density} makes the common calculation visible.

\begin{proposition}[Carlitz's weight as a connector]\label{prop:research-carlitz-bridge}
Let $q_a$ denote the density in \eqref{eq:car-Bernoulli-density},
extended to $-1<a<1$ and interpreted continuously at $a=0$. Then
\begin{equation}\label{eq:research-carlitz-bridge}
 \boxed{\displaystyle
 \Pbb(X_{1+a}\in dy\mid X_2=0)
 =q_a(y)\,dy
 =\frac{C_{1+a}(y)C_{1-a}(y)}{C_2(0)}\,dy.}
\end{equation}
\end{proposition}

\begin{proof}
Only the identification of the density remains. Since $C_2(0)=1/\pi$,
Euler's reflection formula gives
\[
 \left|\Gamma\!\left(\frac{1+a+iy}{2}\right)
        \Gamma\!\left(\frac{1-a+iy}{2}\right)\right|^2
 =\frac{\pi^2}{|\cos(\pi(a+iy)/2)|^2}
 =\frac{2\pi^2}{\cos\pi a+\cosh\pi y}.
\]
Insert the constants in \eqref{eq:Ceta} and use
$\Gamma(1+a)\Gamma(1-a)=\pi a/\sin\pi a$. The result is
\eqref{eq:car-Bernoulli-density} for $a\ne0$, and continuity supplies
$a=0$.
\end{proof}

Carlitz's Euler branch therefore describes an unconditioned Meixner
law, while his Bernoulli/Pasternack branch is obtained by conditioning
the same process on its endpoint. At the midpoint the connector has the
logistic density
\[
 q_0(y)=\frac{\pi}{4\cosh^2(\pi y/2)};
\]
at time $3/2$ it has the hyperbolic-secant density
$q_{1/2}(y)=1/\cosh\pi y$. The symmetry $q_a=q_{-a}$ expresses time
reversal of the zero-endpoint connector, together with its spatial symmetry.
Finally, \eqref{eq:research-bridge-moments} gives
\[
 \operatorname{Var}_{q_a}(Y)=\frac{1-a^2}{3},
\]
exactly the first coefficient of the monic Bernoulli/Pasternack
recurrence. A conditional probability has explained both the family
and its simplest recurrence coefficient.

Carlitz explicitly says that the higher-order Euler construction
was suggested by Stieltjes' continued fraction for the Laplace transform of
$\operatorname{sech}^{\eta}t$, citing Wall's book
\cite[\S9, (9.2)]{carlitz1959bernoulli}. We will derive that
fraction in Section~\ref{sec:stieltjes}, rather than treat it as a recognition
trick. The general relation between recurrences, moment fractions, and Hankel
determinants is established in Section~\ref{sec:fractions}. Carlitz's
contribution is especially illuminating here: familiar numbers lead to a
functional, the functional finds its polynomials, and a suitable realization
of the functional reveals the measure.


\section{Paths, permutations, and the formal fraction}\label{sec:fractions}

The recurrence contains a small combinatorial machine. Write
\begin{equation}\label{eq:b-lambda}
 b_n=(2n+\eta)\delta,\qquad
 \lambda_n=n(n+\eta-1)(1+\delta^2).
\end{equation}
Multiplication by $x$ takes $M_k$ to
$M_{k+1}+b_kM_k+\lambda_kM_{k-1}$. Think of these three terms as
an up step, a level step, and a down step, starting at height $k$.
Their weights are respectively $1$, $b_k$, and $\lambda_k$.
A down step from height zero is forbidden.

After $n$ multiplications, the coefficient $q_{n,k}$ in
\eqref{eq:connection} is the sum of the weights of all paths of length
$n$ from height zero to height $k$. The moment $\mu_n=q_{n,0}$ counts
the returning paths, usually called Motzkin paths. This is the mechanism of Viennot's interpretation of orthogonal polynomials
\cite{viennot1983combinatorial,hetyei2010meixner}. It follows directly
from matrix multiplication, so remains an identity for
weights of both signs. Here all weights are nonnegative when
$\delta\ge0$; for $\delta<0$ the horizontal weights are signed.

\subsection{First return gives the continued fraction}

Let $S_j(z)$ be the formal generating function of excursions based at
height $j$ that never go below $j$. Such an excursion is a sequence of
two possible pieces: a level step of weight $b_jz$, or an arch
consisting of an up step, an excursion at height $j+1$, and a down
step. The arch has weight $\lambda_{j+1}z^2S_{j+1}(z)$. Consequently
\[
 S_j(z)=\frac{1}{1-b_jz-\lambda_{j+1}z^2S_{j+1}(z)}.
\]
Starting at height zero gives
\begin{equation}\label{eq:formal-J}
 \boxed{\displaystyle
 \sum_{n\ge0}\mu_n(\delta,\eta)z^n
 =\cfrac{1}{1-\eta\delta z-
   \cfrac{\eta(1+\delta^2)z^2}{1-(\eta+2)\delta z-
   \cfrac{2(\eta+1)(1+\delta^2)z^2}
         {1-(\eta+4)\delta z-\cdots}}}.}
\end{equation}
The formal meaning is precise: at any specified degree, only finitely
many heights can be reached, so every sufficiently deep truncation
has the same coefficient. This is the path form of the classical
Jacobi fraction correspondence \cite{Flajolet1980,Wall1948}.

It is important to keep this ordinary generating function separate
from the exponential generating function \eqref{eq:H}. In fact,
the ordinary series in \eqref{eq:formal-J} has radius zero. For any
$R>0$, positivity of the density on the whole line gives
\[
 \mu_{2n}\ge R^{2n}\mu_\delta\{|x|>R\}.
\]
Taking $2n$th roots and then allowing $R$ to increase shows
$\mu_{2n}^{1/(2n)}\to\infty$. The continued fraction nevertheless
has a natural analytic meaning after the spectral variable is placed
in its denominators. We will prove that meaning in the next section.

The same recurrence also evaluates the Hankel determinants of our
moment sequence. Let
\[
 \Delta_N=\det[\mu_{j+k}]_{j,k=0}^N.
\]
Changing from the monomials to monic orthogonal polynomials is a
triangular change of basis with determinant one. It diagonalizes this
Gram matrix, so \eqref{eq:orthogonality} gives
\begin{equation}\label{eq:hankel}
 \Delta_N
 =\prod_{j=0}^N j!(\eta)_j(1+\delta^2)^j
 =(1+\delta^2)^{N(N+1)/2}
       \prod_{j=0}^N j!(\eta)_j.
\end{equation}
The product evaluations are another form of the norm calculation,
rather than a separate miracle of determinants.

\subsection{Why powers of the secant count alternating permutations}

At $\delta=0$, horizontal steps disappear. The odd moments vanish,
and the even moments have generating function $\sec^\eta t$.
There is a particularly attractive interpretation of the parameter
$\eta$. Let $\mathcal Z_{2n}$ be the permutations satisfying
\[
 \sigma_1<\sigma_2>\sigma_3<\cdots<\sigma_{2n},
\]
and let $\operatorname{lrmin}(\sigma)$ count their left-to-right
minima, including the first entry. Then
\begin{equation}\label{eq:alternating-moments}
 \mu_{2n}(0,\eta)
 =\sum_{\sigma\in\mathcal Z_{2n}}
       \eta^{\operatorname{lrmin}(\sigma)}.
\end{equation}
This refinement is due to Viennot; Hetyei uses it to connect the
Meixner moments with representation theory
\cite[\S2]{hetyei2010meixner}.

Here is a generating-function proof that also explains the statistic.
Let $E_{2n,k}$ count the permutations with $k$ such minima, and let
$T_{2j+1}$ be the number of alternating permutations of odd length
$2j+1$. Their exponential generating function is $\tan t$.
For completeness, the usual decomposition of an odd alternating
permutation at its largest entry gives
$T'(t)=1+T(t)^2$, $T(0)=0$, and hence $T(t)=\tan t$.

In a permutation from $\mathcal Z_{2n}$, the entry $1$ is in an odd
position, say $2m+1$. The left block has even length $2m$ and carries
$k-1$ left-to-right minima. The entry $1$ is the last new minimum;
the right block has odd length and contributes no further minima.
Its alternating pattern is reversed, but complementation of its
labels gives the same tangent-number count.
Choosing the labels in the left block gives
\[
 E_{2n,k}=\sum_{m=0}^{n-1}\binom{2n-1}{2m}
             E_{2m,k-1}T_{2n-2m-1}.
\]
With
$\Phi_k(t)=\sum_{n\ge0}E_{2n,k}t^{2n}/(2n)!$, this becomes
\[
 \Phi_0(t)=1,\qquad
 \Phi_k'(t)=\Phi_{k-1}(t)\tan t,\qquad \Phi_k(0)=0\quad(k>0).
\]
Integration yields
\[
 \Phi_k(t)=\frac{(-\log\cos t)^k}{k!},\qquad
 \sum_{k\ge0}\eta^k\Phi_k(t)=\sec^\eta t,
\]
proving \eqref{eq:alternating-moments}. For example,
\[
 \mu_2=\eta,\qquad \mu_4=3\eta^2+2\eta,\qquad
 \mu_6=15\eta^3+30\eta^2+16\eta.
\]
At $\eta=1$ these are the secant numbers. At $\eta=2$ they are the
odd tangent numbers, since $(\tan t)'=\sec^2t$. The latter
specialization will bring the even values of the zeta function into
the picture.


\section{What the continued fraction converges to}\label{sec:analytic-fractions}

The analytic companion of the formal moment series is its Cauchy
transform
\begin{equation}\label{eq:Cauchy-x}
 C_\delta(\zeta)=\int_{\R}\frac{K_\eta(x,\delta)}{\zeta-x}\,dx,
 \qquad \zeta\in\C\setminus\R.
\end{equation}
With this sign convention, $C_\delta$ has negative imaginary part
in the upper half-plane. The Nevanlinna convention, which we will
use for the parameter fraction, has the opposite sign.

\subsection{Finite matrices give an infinite limit}

\begin{theorem}\label{thm:Jacobi-fraction}
For $\eta>0$ and real $\delta$, the finite convergents of
\begin{equation}\label{eq:analytic-J}
 C_\delta(\zeta)=
 \cfrac1{\zeta-\eta\delta-
   \cfrac{\eta(1+\delta^2)}{\zeta-(\eta+2)\delta-
   \cfrac{2(\eta+1)(1+\delta^2)}
   {\zeta-(\eta+4)\delta-\cdots}}}
\end{equation}
converge locally uniformly on $\C\setminus\R$ to the integral
\eqref{eq:Cauchy-x}.
\end{theorem}

\begin{proof}
Let $J_N$ be the real symmetric $N\times N$ matrix with diagonal
$b_0,\ldots,b_{N-1}$ and off-diagonal entries
$\sqrt{\lambda_1},\ldots,\sqrt{\lambda_{N-1}}$.
Successive Schur complements give
\begin{equation}\label{eq:finite-J}
 C_{\delta,N}(\zeta)
 =\langle e_0,(\zeta I-J_N)^{-1}e_0\rangle
 =\cfrac1{\zeta-b_0-
   \cfrac{\lambda_1}{\zeta-b_1-
   \cfrac{\lambda_2}{\ddots-
   \cfrac{\lambda_{N-1}}{\zeta-b_{N-1}}}}}.
\end{equation}
The finite-dimensional spectral theorem makes this the Cauchy
transform of a probability measure $\mu_{\delta,N}$ supported on
the eigenvalues of $J_N$.

Its first $2N$ moments agree with those of $\mu_\delta$:
\begin{equation}\label{eq:quadrature-matching}
 \int x^j\,d\mu_{\delta,N}(x)=\mu_j,
 \qquad 0\le j\le2N-1.
\end{equation}
Indeed, the matrix entry $\langle e_0,J_N^je_0\rangle$ is a sum
over walks. A walk that reaches the missing height $N$ and returns
to zero has length at least $2N$. Before that length, truncation
changes nothing. This is Gaussian quadrature seen through the paths.

The common second moment makes the measures tight for $N\ge2$.
Every weakly convergent subsequence has a limit with all the moments
$\mu_j$: for each fixed $j$, use a common higher even moment in
\eqref{eq:quadrature-matching} to control the tails uniformly.
The determinacy proved in Section~\ref{sec:measure} forces this
limit to be $\mu_\delta$. Thus the whole sequence converges weakly
to $\mu_\delta$.

The function $x\mapsto(\zeta-x)^{-1}$ is bounded and continuous
off the real axis, so the Cauchy transforms converge pointwise.
On a compact set at distance $\varepsilon>0$ from $\R$, they and
their derivatives are bounded by $\varepsilon^{-1}$ and
$\varepsilon^{-2}$. A finite-net argument gives local uniform
convergence.
\end{proof}

The denominator of \eqref{eq:finite-J} is
$\det(\zeta I-J_N)=M_N(\zeta;\delta,\eta)$.
For later comparison, define the numerator polynomials by
\begin{equation}\label{eq:associated-polynomials}
 \mathcal N_0=0,\qquad \mathcal N_1=1,\qquad
 \mathcal N_{n+1}=(\zeta-b_n)\mathcal N_n
                         -\lambda_n\mathcal N_{n-1}\quad(n\ge1).
\end{equation}
Then $C_{\delta,N}=\mathcal N_N/M_N$.
These associated polynomials are different from the inverse Sheffer
polynomials of Section~\ref{sec:generating}.

Although the ordinary moment series diverges, its coefficients still
describe the transform at infinity:
\begin{equation}\label{eq:Cauchy-asymptotic}
 C_\delta(\zeta)\sim\sum_{n\ge0}\frac{\mu_n}{\zeta^{n+1}}
\end{equation}
in every closed angular sector separated from the real axis.
To prove any finite order of this expansion, use the finite geometric
identity for $(\zeta-x)^{-1}$ and bound the remainder by an absolute
moment. Thus the formal and analytic interpretations fit together
without treating a divergent series as a convergent one.

\subsection{A hypergeometric evaluation, with the half-plane specified}

The fraction can also be summed explicitly. We use
${}_2F_1(a,b;c;q)=\sum_{j\ge0}(a)_j(b)_jq^j/((c)_jj!)$ near
$q=0$, with its usual continuation cut along $[1,\infty)$.

\begin{proposition}\label{prop:Cauchy-hypergeom}
For $\Im\zeta<0$,
\begin{equation}\label{eq:Cauchy-lower}
 C_\delta(\zeta)=\frac1{\zeta-i\eta}
 {}_2F_1\!\left(1,\eta;
       1+\frac{\eta+i\zeta}{2};\frac{1+i\delta}{2}\right).
\end{equation}
For $\Im\zeta>0$, its conjugate form is
\begin{equation}\label{eq:Cauchy-upper}
 C_\delta(\zeta)=\frac1{\zeta+i\eta}
 {}_2F_1\!\left(1,\eta;
       1+\frac{\eta-i\zeta}{2};\frac{1-i\delta}{2}\right).
\end{equation}
\end{proposition}

\begin{proof}
For $\Im\zeta<0$,
\[
 \frac1{\zeta-x}=i\int_0^\infty e^{-i(\zeta-x)u}\,du.
\]
Absolute Fubini is justified by the integrable majorant
$e^{(\Im\zeta)u}$. The characteristic function obtained from
\eqref{eq:measure-mgf} therefore gives
\begin{equation}\label{eq:Cauchy-Laplace}
 C_\delta(\zeta)
 =i\int_0^\infty e^{-i\zeta u}
           (\cosh u-i\delta\sinh u)^{-\eta}\,du.
\end{equation}
Put $r=e^{-2u}$, $a=(\eta+i\zeta)/2$, and
$q=(1+i\delta)/2$. Then $\Re a>0$ and
\begin{align*}
 C_\delta(\zeta)
 &=\frac i2\int_0^1 r^{a-1}[(1-q)+qr]^{-\eta}\,dr\\
 &=\frac{i}{2a}(1-q)^{-\eta}
    {}_2F_1\!\left(\eta,a;a+1;-\frac q{1-q}\right)\\
 &=\frac{i}{2a}\,{}_2F_1(\eta,1;a+1;q).
\end{align*}
The second line is Euler's beta-integral representation; the third
is Pfaff's transformation. Both follow with the branches continued
from small $q$ to $\Re q=1/2$, where no cut is crossed
\cite[\S\S15.6, 15.8]{DLMF}.
Since $i/(2a)=1/(\zeta-i\eta)$, this proves
\eqref{eq:Cauchy-lower}. Conjugating the real density proves
\eqref{eq:Cauchy-upper}.
\end{proof}

The proof gives a summation of the Jacobi fraction by identifying its
analytic function first. It also prepares a useful surprise: the
same hypergeometric expression will appear when $\delta$, rather
than $x$, becomes the spectral variable.

\subsection{Stieltjes' fraction and Carlitz's starting point}\label{sec:stieltjes}

At $\delta=0$, the fraction takes a familiar form. For $\Re s>0$,
let
\[
 L_\eta(s)=\int_0^\infty e^{-su}\operatorname{sech}^{\eta}u\,du.
\]
Since $C_0(-is)=iL_\eta(s)$, Theorem~\ref{thm:Jacobi-fraction}
gives the convergent continued fraction
\begin{equation}\label{eq:Stieltjes-sech}
 \boxed{\displaystyle
 L_\eta(s)=\cfrac1{s+
   \cfrac{\eta}{s+
   \cfrac{2(\eta+1)}{s+
   \cfrac{3(\eta+2)}{s+\cdots}}}},
 \qquad \Re s>0.}
\end{equation}
The convergence is locally uniform in this half-plane. This is the
Stieltjes fraction that Carlitz cites as the impetus for his
higher-order Euler construction \cite[\S9]{carlitz1959bernoulli}.

Its coefficients can be derived directly by integration by parts.
For real $s>0$, define
\[
 A_n(s)=\int_0^\infty e^{-su}\tanh^n u\,
                         \operatorname{sech}^{\eta}u\,du.
\]
Differentiating $\tanh^n u\operatorname{sech}^{\eta}u$ gives
\[
 \frac d{du}\bigl(\tanh^n u\operatorname{sech}^{\eta}u\bigr)
 =n\tanh^{n-1}u\operatorname{sech}^{\eta}u
 -(n+\eta)\tanh^{n+1}u\operatorname{sech}^{\eta}u.
\]
The endpoint terms yield
\begin{equation}\label{eq:A-integral-recurrence}
 1=sA_0+\eta A_1,\qquad
 nA_{n-1}=sA_n+(n+\eta)A_{n+1}\quad(n\ge1).
\end{equation}
Since $A_n>0$, their ratios may be substituted successively:
\[
 \frac{A_n}{A_{n-1}}
 =\frac{n}{s+(n+\eta)A_{n+1}/A_n}.
\]
Together with the first identity, this produces exactly
\eqref{eq:Stieltjes-sech}. The preceding Jacobi argument supplies
the convergence that a succession of formal substitutions alone
would not establish. The two derivations have different jobs: the
integrals explain the coefficients, and the finite matrices prove
the infinite limit.


\section{A tangent-number encore: the even zeta values}\label{sec:zeta}

The specialization $\eta=2$, $\delta=0$ has a particularly simple
density. The reflection formula for the gamma function, or the
convolution in Carlitz's construction, gives
\begin{equation}\label{eq:tangent-density}
 K_2(x,0)=\frac{x}{2\sinh(\pi x/2)},
\end{equation}
with its continuous value $1/\pi$ at $x=0$. The convolution is
worth seeing once. If $p_1(x)=1/(2\cosh(\pi x/2))$ and
$v=\pi x/2$, then
\[
 (p_1*p_1)(x)=\frac1{2\pi}
       \int_{\R}\frac{du}{\cosh u\cosh(v-u)}.
\]
The identity
\[
 \frac1{\cosh u\cosh(v-u)}
 =\frac{\tanh u+\tanh(v-u)}{\sinh v}
\]
reduces the integral to $2v/\sinh v$ by evaluating
$\log\cosh u-\log\cosh(v-u)$ at the endpoints. This proves
\eqref{eq:tangent-density} directly from Carlitz's two summands.

Write
$\tan t=\sum_{m\ge0}T_{2m+1}t^{2m+1}/(2m+1)!$.
Since $\sec^2t=(\tan t)'$, our moment calculation says
\[
 T_{2m+1}=\int_{\R}x^{2m}K_2(x,0)\,dx.
\]
For $x>0$, expand
\[
 \frac1{2\sinh(\pi x/2)}
 =\sum_{j\ge0}e^{-(j+1/2)\pi x}.
\]
Every term is nonnegative, so termwise integration follows from
monotone convergence. The elementary gamma integral then gives
\begin{align}\label{eq:tangent-zeta}
 T_{2m+1}
 &=2\sum_{j\ge0}\int_0^\infty
                 x^{2m+1}e^{-(j+1/2)\pi x}\,dx\notag\\
 &=\frac{2(2m+1)!}{\pi^{2m+2}}
       (2^{2m+2}-1)\zeta(2m+2).
\end{align}
The odd integers in the geometric expansion explain the factor
$2^{2m+2}-1$.

To recover Euler's usual form, define the Bernoulli numbers by
\[
 \frac{z}{e^z-1}=\sum_{n\ge0}B_n\frac{z^n}{n!}.
\]
Substituting $z=2it$ and using
$\tan t=\cot t-2\cot(2t)$ gives, by coefficient comparison,
\begin{equation}\label{eq:tangent-Bernoulli}
 T_{2r-1}=(-1)^{r-1}
       \frac{2^{2r}(2^{2r}-1)}{2r}B_{2r},\qquad r\ge1.
\end{equation}
Combining \eqref{eq:tangent-zeta} and
\eqref{eq:tangent-Bernoulli} yields
\begin{equation}\label{eq:Euler-zeta}
 \boxed{\displaystyle
 \zeta(2r)=(-1)^{r+1}
             \frac{(2\pi)^{2r}B_{2r}}{2(2r)!}.}
\end{equation}
The measure was obtained from a beta integral, so this argument does
not use the zeta evaluation it is proving. It places a classical
formula at the meeting point of a positive measure, alternating
permutations, and Carlitz's Euler umbra. For other elementary routes
to the same values, see \cite{Holtz2024Zeta}; for the arithmetic of
tangent and secant numbers, see \cite{BrentHarvey2013}.

\subsection{The same density measures the jumps}

The tangent density also has a role in the Meixner process of
Section~\ref{sec:meixner-time}. Put $\theta=\arctan\delta$ and
write its cumulant exponent per unit time as
\[
 \Psi_\delta(u)=-\log(\cos u-\delta\sin u),
 \qquad \log\E e^{uX_t}=t\Psi_\delta(u).
\]
The positive measure
\begin{equation}\label{eq:research-levy-density}
 \nu_\delta(dy)
 =\frac{e^{\theta y}}{2y\sinh(\pi y/2)}\,dy
 =(1+\delta^2)\frac{K_2(y,\delta)}{y^2}\,dy
\end{equation}
is its \emph{L\'evy measure}: it records the intensity of the jumps.
This classical formula \cite[\S3.1]{MazzolaMuliere2011} follows
directly from the density already in hand. Indeed,
\[
 y^2\nu_\delta(dy)=e^{\theta y}C_2(y)\,dy,
 \qquad
 \int_{\R}y^2e^{uy}\,\nu_\delta(dy)
       =\sec^2(\theta+u)=\Psi_\delta''(u).
\]
The measure satisfies $\int\min(1,y^2)\,\nu_\delta(dy)<\infty$:
its density is asymptotic to $1/(\pi y^2)$ at zero and has
exponentially decreasing tails. Integrating twice and fixing the value
and derivative at zero therefore gives
\begin{equation}\label{eq:research-levy-exponent}
 \Psi_\delta(u)
 =\delta u+\int_{\R}(e^{uy}-1-uy)\,\nu_\delta(dy).
\end{equation}
The compensation cancels the singularity at zero, and the integral
converges throughout $|\theta+u|<\pi/2$. This is the
L\'evy--Khintchine representation; no Gaussian term is present.

If $\kappa_n(\delta,\eta)$ denotes the $n$th cumulant of our law,
then $\log H(u)=\eta\Psi_\delta(u)$, and differentiation gives
\begin{equation}\label{eq:research-cumulant-moment}
 \boxed{\displaystyle
 \kappa_{n+2}(\delta,\eta)
 =\eta(1+\delta^2)\mu_n(\delta,2),\qquad n\ge0.}
\end{equation}
Thus the shape-two moments are the higher cumulants of the entire
Meixner convolution family. At zero tilt they are the tangent numbers
just used to evaluate the even zeta values. The same small density
connects moments, cumulants, and jumps.


\section{Turning the recurrence sideways}\label{sec:parameter-duality}

So far $x$ has been the variable and $\delta$ a parameter. Nothing in
the recurrence prevents us from reversing those roles. The question
is what kind of orthogonality survives the reversal. The factor
$1+\delta^2$ tells us in advance that the answer will involve
rational functions.

\subsection{The parameter recurrence}

Fix $x\in\R$ and define
\begin{equation}\label{eq:RII-normalization}
 P_n(z)=\frac{(-1)^nM_n(x;z,\eta)}{2^n(\eta/2)_n}.
\end{equation}
The recurrence becomes
\begin{equation}\label{eq:RII-recurrence}
 P_{n+1}(z)=(z-r_n)P_n(z)-c_n(1+z^2)P_{n-1}(z),
\end{equation}
where
\begin{equation}\label{eq:RII-coefficients}
 r_n=\frac{x}{2n+\eta},\qquad
 c_n=\frac{n(n+\eta-1)}{(2n+\eta-2)(2n+\eta)}\quad(n\ge1).
\end{equation}
Here $P_0=1$ and $P_1=z-r_0$. These polynomials are not generally
monic in $z$: the quadratic factor also contributes to the leading
coefficient. A recurrence with a quadratic factor multiplying the
preceding polynomial is of type $R_{\mathrm{II}}$; a linear factor
would give type $R_{\mathrm I}$ \cite{IsmailMasson1995}.

\subsection{Why Wall's fraction is the right question}\label{sec:Wall-introduction}

The Jacobi fraction in Section~\ref{sec:analytic-fractions} encoded
moments through an expansion at infinity. The factor $1+z^2$ suggests
a different expansion point: it vanishes at $i$ and $-i$. In the
rational orthogonality these same points will be the fixed poles
of the test functions. They are natural when the real line is
viewed as a circle. The
Cayley map
\begin{equation}\label{eq:Wall-Cayley}
 \zeta=\frac{i-z}{i+z},\qquad
 z=i\frac{1-\zeta}{1+\zeta}
\end{equation}
sends the upper half-plane to the unit disc, $i$ to its center, and
the real line to its boundary, with infinity represented by $-1$.
Taylor interpolation at $i$ is therefore the half-plane counterpart
of Taylor interpolation at the center of the disc. The conjugate
point $-i$ enters when the approximants have real coefficients.

Here is the measure-theoretic part of the picture. A
\emph{Nevanlinna function}, also called a Herglotz function, is an
analytic function $h$ on $\C_+$ with positive imaginary part. Such
a function has the representation
\begin{equation}\label{eq:general-Nevanlinna}
 h(z)=az+b+\int_{\R}\frac{1+tz}{t-z}\,d\sigma(t),
 \qquad a\ge0,\quad b\in\R,\quad \sigma\ge0,\quad\sigma(\R)<\infty.
\end{equation}
The measure here is finite. In the usual subtracted Cauchy
representation the measure is $d\rho(t)=(1+t^2)\,d\sigma(t)$,
which need not have finite mass, and the integral becomes
\[
 \int_{\R}\left(\frac1{t-z}-\frac{t}{1+t^2}\right)d\rho(t).
\]
The subtraction ensures convergence when the two terms cannot be
integrated separately. Positivity is visible in the original kernel:
\begin{equation}\label{eq:Nevanlinna-kernel-imaginary}
 \Im\frac{1+tz}{t-z}
 =\frac{(1+t^2)\Im z}{|t-z|^2}.
\end{equation}

For completeness, the appearance of a measure can be seen directly
on the disc. Put $H(\zeta)=-ih(z)$ with $z$ as in
\eqref{eq:Wall-Cayley}. Then $u=\Re H$ is positive and harmonic.
The measures $u(re^{i\theta})\,d\theta/(2\pi)$ have the common
mass $u(0)$. A weakly convergent subsequence as $r\uparrow1$ gives
a positive boundary measure $\mu$. The Poisson formula on each
smaller disc, followed by this limit, represents $u$ as its Poisson
integral. Adding the harmonic conjugate gives
\[
 H(\zeta)=i\Im H(0)+
          \int_{|\xi|=1}\frac{\xi+\zeta}{\xi-\zeta}\,d\mu(\xi).
\]
For $\xi_t=(i-t)/(i+t)$, a calculation gives
\[
 i\frac{\xi_t+\zeta}{\xi_t-\zeta}
       =\frac{1+tz}{t-z}.
\]
The atom at $\xi=-1$ contributes $az$, and the rest of $\mu$
pushes forward to $\sigma$ on $\R$. This proves
\eqref{eq:general-Nevanlinna}; the real constant is
$b=\Re h(i)$. Uniqueness follows from the Fourier coefficients of
the boundary measure, or from the Poisson formula itself.

At $z=i$, the integral kernel is exactly $i$. Thus
\begin{equation}\label{eq:Nevanlinna-normalization}
 h(i)=i\quad\Longleftrightarrow\quad
 b=0,\qquad a+\sigma(\R)=1.
\end{equation}
A normalized Nevanlinna function therefore represents a probability
on the real line together with a possible atom at infinity. It
represents a probability on $\R$ itself when $a=0$. The coefficient
of the atom is
\[
 a=\lim_{y\to\infty}\frac{\Im h(iy)}{y},
\]
as follows from \eqref{eq:general-Nevanlinna} by dominated
convergence. This small qualification will matter when we compare
normalizations.

Wall found a continued-fraction description of this correspondence
\cite{Wall1944,Wall1946,Wall1948}. In the disc, a function $H$ with
positive real part and $H(0)=1$ can be written
\[
 H(\zeta)=\frac{1+\zeta f(\zeta)}{1-\zeta f(\zeta)},
 \qquad |f(\zeta)|\le1.
\]
Indeed, $(H-1)/(H+1)$ maps the disc into itself and vanishes at zero,
so Schwarz's lemma supplies the factor $\zeta$. The same argument
can be repeated: if $\gamma=f(0)$ and $|\gamma|<1$, then
\[
 f_1(\zeta)=\frac{f(\zeta)-\gamma}
                 {\zeta[1-\overline\gamma f(\zeta)]}
\]
again maps the disc into its closure. Iterating this step is the
Schur algorithm; if $|\gamma|=1$, the maximum principle makes
$f$ constant and the procedure terminates. Wall reorganized the
resulting fractions. After the Cayley change of variable, his
canonical half-plane form is
\begin{equation}\label{eq:general-Wall-model}
 \Phi(z)=-\cfrac1{z-
   \cfrac{g_1(1+z^2)}{z-r_1-
   \cfrac{(1-g_1)g_2(1+z^2)}{z-r_2-
   \cfrac{(1-g_2)g_3(1+z^2)}{z-r_3-\cdots}}}}.
\end{equation}
For real $r_n$ and $0<g_n<1$, this fraction converges locally
uniformly on $\C_+$ to a Nevanlinna function with $\Phi(i)=i$.
This is the form of Wall's theorem presented by Derevyagin
\cite[Theorem~3.1, formula~(3.12), and Corollary~3.2]{Derevyagin2017},
where the associated integral representation is also displayed in
Theorem~5.3. The first denominator $z$ records the normalization
at $i$; the factors $1+z^2$ record the interpolation points.

Our coefficients have exactly Wall's chain form:
\begin{equation}\label{eq:chain-sequence}
 c_n=(1-g_{n-1})g_n,\qquad
 g_n=\frac{n}{2n+\eta},\qquad g_0=0.
\end{equation}
In particular $c_1=g_1$, and $0<g_n<1$ for every $n\ge1$ and
every $\eta>0$. The only initial difference is the shift $r_0$
in the first denominator. We shall first identify the transform
with that natural shift, deriving its measure and convergence
directly. Afterwards a single fractional linear transformation
will put it into \eqref{eq:general-Wall-model}. The general theory
has now told us what to look for; the following calculations prove
the needed instance.

\subsection{The positive pencil}

The two spectral directions can already be seen in finite matrices.
Index rows by $k=0,\ldots,N-1$, and put
$a_k=\sqrt{(k+1)(k+\eta)}$. Let $A_N$ have diagonal $2k+\eta$
and upper and lower entries $a_k$. Let $B_N$ have zero diagonal,
upper entries $-ia_k$, and lower entries $ia_k$. Expansion of the
tridiagonal determinant gives
\begin{equation}\label{eq:two-variable-determinant}
 \boxed{\displaystyle
 M_N(x;z,\eta)=\det(xI_N-B_N-zA_N).}
\end{equation}
Moreover $A_N=L_NL_N^*$, where $L_N$ is lower bidiagonal with
diagonal $\sqrt{k+\eta}$ and subdiagonal $\sqrt{k}$. Thus $A_N$
is positive definite, while $B_N$ is Hermitian.

For real $z=\delta$, the matrix $B_N+\delta A_N$ is Hermitian.
A diagonal unitary change of basis makes its off-diagonal entries
$a_k\sqrt{1+\delta^2}$, recovering the Jacobi matrix of
Section~\ref{sec:analytic-fractions}. For fixed real $x$, the same
determinant describes the generalized eigenvalue problem
\begin{equation}\label{eq:generalized-eigenvalue}
 (xI_N-B_N)v=zA_Nv.
\end{equation}
Its eigenvalues are real because $A_N>0$: conjugating by
$A_N^{-1/2}$ gives an ordinary Hermitian problem. In particular,
the zeros of $M_N$ in the parameter are real, just as its zeros in
the argument are real. Related unit-circle and pencil constructions
are developed in \cite{Derevyagin2017,IsmailSriRanga2019}.

\subsection{A Rodrigues formula in the parameter}

Keep $x$ fixed and write
\begin{equation}\label{eq:dual-weight}
 s(t)=1+t^2,\qquad
 w(t)=s(t)^{-\eta/2}e^{x\arctan t},\qquad
 K_\eta(x,t)=C_\eta(x)w(t).
\end{equation}
The logarithmic derivative is remarkably simple:
\begin{equation}\label{eq:dual-Pearson}
 s(t)w'(t)=(x-\eta t)w(t).
\end{equation}
The generating function gives a Rodrigues formula without any
guesswork. Taylor-expand $w(t+us(t))/w(t)$. The elementary identities
\[
 \frac{1+(t+us(t))^2}{1+t^2}
       =1+2tu+(1+t^2)u^2,
 \qquad
 \arctan(t+us(t))-\arctan t
       =\arctan\frac{u}{1+tu}
\]
hold near $u=0$. Consequently this Taylor series is exactly
$G(x;u)$ with $\delta=t$, and
\begin{equation}\label{eq:dual-Rodrigues}
 \boxed{\displaystyle
 M_n(x;t,\eta)=\frac{s(t)^n}{w(t)}\frac{d^nw(t)}{dt^n}.}
\end{equation}
These are the complementary Romanovski--Routh polynomials in the
normalization of our recurrence
\cite{MartinezFinkelshteinEtAl2019}. The identification explains
the family; the calculation proves the formula in the conventions
used here.

The Rodrigues formula also gives a lowering operation in the
parameter. Differentiate it, use \eqref{eq:dual-Pearson}, and
recognize the next Rodrigues expression:
\[
 \partial_tM_n
 =\frac{M_{n+1}-[x-(\eta+2n)t]M_n}{1+t^2}.
\]
The recurrence therefore reduces the derivative to
\begin{equation}\label{eq:delta-lowering-expanded}
 \partial_\delta M_n=-n(n+\eta-1)M_{n-1}.
\end{equation}
Thus $(-1)^nM_n/(\eta)_n$ is an Appell sequence in the parameter.
Differentiate once more and use the recurrence at index $n-1$ to
obtain
\begin{equation}\label{eq:delta-ODE-expanded}
 (1+\delta^2)\partial_\delta^2M_n
 +[x-(\eta+2n-2)\delta]\partial_\delta M_n
 +n(n+\eta-1)M_n=0.
\end{equation}
For $n\ge2$, if $c=n(n+\eta-1)$ and
$e=(n-1)(n+\eta-2)$, the left side is
\[
 c\{(1+\delta^2)eM_{n-2}
   -[x-(\eta+2n-2)\delta]M_{n-1}+M_n\}=0.
\]
The cases $n=0,1$ follow directly. This also derives the
complementary Romanovski--Routh differential equation in our
normalization.

Integration by parts now gives rational orthogonality:
\begin{equation}\label{eq:rational-orthogonality}
 \int_{\R}\frac{t^kM_n(x;t,\eta)}{(1+t^2)^n}
                K_\eta(x,t)\,dt=0,
 \qquad 0\le k<n,\quad\eta>0.
\end{equation}
Indeed, the integral is $C_\eta(x)\int t^kw^{(n)}(t)\,dt$.
At either endpoint $w^{(j)}(t)=O(|t|^{-\eta-j})$.
Every boundary term after integration by parts is therefore
$O(|t|^{k-n+1-\eta})$, which tends to zero for $k<n$.
After $n$ integrations the derivative of $t^k$ vanishes.
The integrand itself is absolutely integrable, since its growth is
$O(|t|^{k-n-\eta})$.

This is the promised second orthogonality. It is not ordinary
polynomial orthogonality for a fixed measure with all moments:
the denominator grows with the degree. That denominator is exactly
what makes the heavy-tailed parameter measure useful at every degree.

\subsection{Two normalization integrals}

We will need the mass of the parameter weight and its mass after
division by $1+t^2$. Both follow from one Fourier calculation.
For $\Re a>1/2$, the gamma integral gives the Fourier transform
\[
 \int_0^\infty \frac{y^{a-1}e^{-y}}{\Gamma(a)}e^{-ity}\,dy
       =(1+it)^{-a}.
\]
Plancherel's identity and another gamma integral therefore give
\begin{equation}\label{eq:dual-beta}
 \int_{\R}|(1+it)^{-a}|^2\,dt
 =\frac{2\pi}{|\Gamma(a)|^2}
       \int_0^\infty y^{2\Re a-2}e^{-2y}\,dy
 =\frac{\pi\Gamma(2\Re a-1)}
        {2^{2\Re a-2}|\Gamma(a)|^2}.
\end{equation}
This is a beta integral in Fourier form
\cite[\S5.12]{DLMF}. Choose $a=(\eta+ix)/2$ for the first
identity below, and $a=1+(\eta+ix)/2$ for the second. Multiplying
by $C_\eta(x)$ and cancelling gamma factors yields
\begin{align}
 \int_{\R}K_\eta(x,t)\,dt
     &=\frac1{\eta-1}, &&\eta>1,
       \label{eq:dual-total-mass}\\
 \int_{\R}\frac{K_\eta(x,t)}{1+t^2}\,dt
     &=\frac{\eta}{x^2+\eta^2}, &&\eta>0.
       \label{eq:dual-weighted-mass}
\end{align}
The first integral is infinite for $0<\eta\le1$, since
$w(t)$ has tails proportional to $|t|^{-\eta}$.
When $\eta>1$, we may introduce the probability measure
\begin{equation}\label{eq:nu-dual}
 d\nu_x(t)=(\eta-1)K_\eta(x,t)\,dt.
\end{equation}
The preceding integration by parts then also includes the endpoint
$k=n$:
\begin{equation}\label{eq:rational-triangular}
 \int_{\R}\frac{t^kM_n(x;t,\eta)}{(1+t^2)^n}\,d\nu_x(t)
       =(-1)^n n!\,\delta_{kn},\qquad 0\le k\le n.
\end{equation}

\subsection{The parameter fraction and its exact measure}

For $z$ in the upper half-plane $\C_+$, define
\begin{equation}\label{eq:F-integral-expanded}
 F_x(z)=\eta\int_{\R}\frac{K_\eta(x,t)}{t-z}\,dt.
\end{equation}
This converges for every $\eta>0$, even when the measure in its
numerator has infinite total mass: the integrand is
$O(|t|^{-\eta-1})$. Moreover
\begin{equation}\label{eq:F-imaginary}
 \Im F_x(z)=\eta\Im z
       \int_{\R}\frac{K_\eta(x,t)}{|t-z|^2}\,dt>0.
\end{equation}
Thus $F_x$ is a Nevanlinna (or Herglotz) function.

\begin{theorem}\label{thm:Wall-fraction-expanded}
For every real $x$ and $\eta>0$,
\begin{equation}\label{eq:Wall-fraction-expanded}
 \boxed{\displaystyle
 F_x(z)=-\cfrac1{z-r_0-
       \cfrac{c_1(1+z^2)}{z-r_1-
       \cfrac{c_2(1+z^2)}{z-r_2-\cdots}}},
       \qquad z\in\C_+.}
\end{equation}
The finite convergents converge locally uniformly on $\C_+$ to
the integral \eqref{eq:F-integral-expanded}.
\end{theorem}

\begin{proof}
We first find an equation that characterizes the integral. Put
$I_j=\int t^jw(t)/(1+t^2)\,dt$ for $j=0,1$ and
$\kappa=\eta C_\eta(x)$. Integrating \eqref{eq:dual-Pearson}
after division by $s(t)$ gives $xI_0=\eta I_1$, since $w$ vanishes
at both endpoints. Formula \eqref{eq:dual-weighted-mass} gives
$\kappa I_0=\eta^2/(x^2+\eta^2)$. In particular,
\begin{equation}\label{eq:F-at-i-expanded}
 F_x(i)=\kappa(I_1+iI_0)=\frac{\eta}{x-i\eta}.
\end{equation}

Integration by parts in $w(t)/(t-z)$ gives
$F_x'(z)=\kappa\int w'(t)/(t-z)\,dt$.
Substituting \eqref{eq:dual-Pearson} and cancelling $t-z$ gives
\begin{align*}
 (1+z^2)F_x'(z)-(x-\eta z)F_x(z)
 &=-\kappa[\eta I_0+xI_1+z(xI_0-\eta I_1)]\\
 &=-\eta.
\end{align*}
All boundary terms vanish for $\eta>0$. Thus
\begin{equation}\label{eq:F-ODE-expanded}
 (1+z^2)F_x'(z)=(x-\eta z)F_x(z)-\eta.
\end{equation}
Differentiate $n$ times, for $n\ge1$:
\[
 (1+z^2)F_x^{(n+1)}
 =[x-(\eta+2n)z]F_x^{(n)}
       -n(n+\eta-1)F_x^{(n-1)}.
\]
It follows that $U_n=(1+z^2)^nF_x^{(n)}$ satisfies the original
Meixner recurrence in $n$. If $N_0=0$, $N_1=1$, and $N_n$
satisfies that recurrence with $\delta=z$ for $n\ge1$, the initial
values $U_0=F_x$, $U_1=M_1F_x-\eta$ give the exact identity
\begin{equation}\label{eq:Wall-exact-error}
 M_n(x;z,\eta)F_x(z)-\eta N_n(z)
       =(1+z^2)^nF_x^{(n)}(z).
\end{equation}
The quotient $F_{x,n}=\eta N_n/M_n$ is the $n$-level fraction
in \eqref{eq:Wall-fraction-expanded}; this follows by the ordinary
continuant recurrence, or by taking the first Schur complement in
\eqref{eq:two-variable-determinant}.

It remains to prove convergence. The same Schur complement gives
\begin{equation}\label{eq:Wall-finite-resolvent}
 F_{x,N}(z)=\eta\langle e_0,
                  (xI_N-B_N-zA_N)^{-1}e_0\rangle.
\end{equation}
If $H=xI_N-B_N-zA_N$, then $\Im H=-(\Im z)A_N<0$,
so $H$ is invertible and
$\Im H^{-1}=(\Im z)H^{-*}A_NH^{-1}>0$.
Thus each of these functions has positive imaginary part in $\C_+$.
At $z=i$ the matrix is lower triangular, and
\[
 F_{x,N}(i)=\frac{\eta}{x-i\eta}=F_x(i).
\]
Herglotz functions with this fixed value form a normal family:
map their common value to zero by a half-plane-to-disc map and
apply Schwarz's lemma, followed by Montel's theorem. The fixed
interior value also prevents a subsequence from tending to infinity.

Finally,
\[
 M_n(x;i,\eta)=\prod_{j=0}^{n-1}[x-i(\eta+2j)]\ne0.
\]
Equation \eqref{eq:Wall-exact-error} shows that $F_{x,n}$ and
$F_x$ have the same first $n$ Taylor coefficients at $i$. Every
locally uniform subsequential limit therefore agrees with $F_x$
to every order at $i$. The identity theorem identifies that limit
with $F_x$ throughout $\C_+$. Since every subsequential limit is
the same, the whole sequence converges.
\end{proof}

This proof makes the role of positivity quite concrete. The recurrence
alone produces rational approximants; the positive pencil makes
them a controlled analytic family; the differential identity selects
their unique limit.

Expansion of \eqref{eq:F-ODE-expanded} at $z=i$ also gives
\begin{equation}\label{eq:F-hypergeometric-expanded}
 F_x(z)=\frac{\eta}{x-i\eta}
 {}_2F_1\!\left(1,\eta;
       1+\frac{\eta+ix}{2};\frac{1+iz}{2}\right).
\end{equation}
For example, with $q=(1+iz)/2$ and
$F_x=\sum_{j\ge0}a_jq^j$, the equation gives
$a_0=\eta/(x-i\eta)$ and
$a_{j+1}/a_j=(\eta+j)/(1+(\eta+ix)/2+j)$.
This is precisely the displayed hypergeometric series.
Analytic continuation covers $\C_+$ because $\Re q<1/2$ there.

The match with \eqref{eq:Cauchy-lower} is now exact: the same finite
continued fraction, viewed in $x$, is an ordinary Jacobi fraction;
viewed in $\delta$, it is $1/\eta$ times this Wall fraction.
The lower-half-plane branch in the first variable leads to the
upper-half-plane Nevanlinna function in the second. The distinction
between the two branches is essential when identifying their
boundary measures.

One can make this analytic connection without merely comparing
formulas. The function
\begin{equation}\label{eq:two-variable-transform}
 \mathcal C(\zeta,z)=i\int_0^\infty e^{-i\zeta u}
                  (\cosh u-iz\sinh u)^{-\eta}\,du
\end{equation}
is jointly holomorphic for $\Im\zeta<\eta$, $\Im z>-1$.
Indeed, the expression in parentheses has positive real part and,
locally uniformly in this domain, its modulus is bounded below
by a positive constant times $\cosh u$. We use its principal power.
Equation \eqref{eq:Cauchy-Laplace} identifies
$\mathcal C(\zeta,\delta)=C_\delta(\zeta)$ for real $\delta$
and $\Im\zeta<0$; the hypergeometric formulas identify
$\eta\mathcal C(x,z)=F_x(z)$ for real $x$ and $z\in\C_+$.
In particular, the compatible boundary values satisfy
\[
 \eta C_\delta(x-i0)=F_x(\delta+i0),\qquad
 \Im C_\delta(x-i0)=\pi K_\eta(x,\delta).
\]
The upper $x$ boundary is the complex conjugate. Thus the same
analytic continuation carries the two spectral readings, with
their boundary branches explicitly specified.


\section{Two measures, or one measure in two variables?}\label{sec:joint-measure}

Both descriptions are useful. One positive kernel gives the two
spectral realizations, and a Cauchy weighting turns them into the
conditional laws of one joint probability. The factor $1+t^2$
determines which parameter measure appears: a Cauchy integral and a
normalized Nevanlinna integral use closely related, but different,
measures.

\subsection{The Cauchy measure and the finite Nevanlinna measure}

The measure in \eqref{eq:F-integral-expanded} is
\begin{equation}\label{eq:rho-expanded}
 d\rho_x(t)=\eta K_\eta(x,t)\,dt.
\end{equation}
Its mass is $\eta/(\eta-1)$ when $\eta>1$ and infinite when
$0<\eta\le1$. Nevertheless the Cauchy integral is absolutely
convergent for every $\eta>0$.

For the Nevanlinna kernel, divide this measure by $1+t^2$:
\begin{equation}\label{eq:sigma-expanded}
 d\sigma_x(t)=\frac{d\rho_x(t)}{1+t^2}
       =\eta\frac{K_\eta(x,t)}{1+t^2}\,dt.
\end{equation}
The identity
\[
 \frac1{t-z}=\frac1{1+t^2}
          \left(\frac{1+tz}{t-z}+t\right)
\]
and the value \eqref{eq:F-at-i-expanded} give
\begin{equation}\label{eq:F-Nevanlinna}
 F_x(z)=\frac{\eta x}{x^2+\eta^2}
       +\int_{\R}\frac{1+tz}{t-z}\,d\sigma_x(t),
 \qquad \sigma_x(\R)=\frac{\eta^2}{x^2+\eta^2}.
\end{equation}
There is no linear term. More directly than an appeal to the
general representation theorem, dominated convergence gives
$F_x(iy)\to0$ as $y\to\infty$: for $y\ge1$,
$|t-iy|^{-1}\le(1+t^2)^{-1/2}$, and
$K_\eta(x,t)(1+t^2)^{-1/2}$ is integrable for every $\eta>0$.

To normalize at $i$, set
\begin{equation}\label{eq:phi-affine}
 \varphi_x(z)=\frac{x^2+\eta^2}{\eta^2}F_x(z)-\frac{x}{\eta}.
\end{equation}
Then $\varphi_x(i)=i$ and
\begin{equation}\label{eq:phi-probability}
 \boxed{\displaystyle
 \varphi_x(z)=\int_{\R}\frac{1+tz}{t-z}\,d\widehat\sigma_x(t),
 \qquad
 d\widehat\sigma_x(t)=\frac{x^2+\eta^2}{\eta}
             \frac{K_\eta(x,t)}{1+t^2}\,dt.}
\end{equation}
The measure $\widehat\sigma_x$ is a probability for all $\eta>0$.
Gamma cancellation in \eqref{eq:Ceta} also gives
\begin{equation}\label{eq:shifted-eta-measure}
 d\widehat\sigma_x(t)=(\eta+1)K_{\eta+2}(x,t)\,dt.
\end{equation}
Thus this normalized probability is the unweighted parameter
probability with the shape parameter shifted by two. For the
original polynomials, its rational orthogonality is equivalently
\begin{equation}\label{eq:weighted-rational-orthogonality}
 \int_{\R}\frac{t^kM_n(x;t,\eta)}{(1+t^2)^{n-1}}
       \,d\widehat\sigma_x(t)=0,
 \qquad n\ge1,\quad 0\le k<n.
\end{equation}

\subsection{The canonical Wall normalization}

We can now return to Wall's canonical fraction
\eqref{eq:general-Wall-model}. Its first denominator is $z$;
the natural Meixner fraction \eqref{eq:Wall-fraction-expanded}
starts with $z-r_0$, where $r_0=x/\eta$. Removing just that
shift means adding $r_0$ to the reciprocal denominator:
\[
 -\frac1{W_x(z)}=-\frac1{F_x(z)}+r_0.
\]
Solving this equation derives the required transformation:
\begin{equation}\label{eq:canonical-Wall}
 W_x(z)=\frac{F_x(z)}{1-(x/\eta)F_x(z)}
       =-\cfrac1{z-
          \cfrac{c_1(1+z^2)}{z-r_1-
          \cfrac{c_2(1+z^2)}{z-r_2-\cdots}}}.
\end{equation}
The same identity holds for every finite convergent, so the
local uniform convergence proved in
Theorem~\ref{thm:Wall-fraction-expanded} also justifies this
infinite fraction. It retains the entire continued-fraction tail.
The denominator in the fractional linear transformation cannot
vanish on $\C_+$,
since $F_x$ has strictly positive imaginary part. Moreover
\[
 \Im W_x(z)=\frac{\Im F_x(z)}{|1-(x/\eta)F_x(z)|^2}>0,
 \qquad W_x(i)=i.
\]
The latter identity follows at once from
$F_x(i)=\eta/(x-i\eta)$. Thus this is precisely the canonical
normalization described by Wall and Derevyagin, for every
$\eta>0$.

There are now two natural ways to normalize at $i$.
The affine transformation \eqref{eq:phi-affine} merely multiplies
the representing measure by a constant, producing the explicit
probability $\widehat\sigma_x$. The fractional linear
transformation \eqref{eq:canonical-Wall} instead preserves the
continued-fraction tail. These operations generally give different
functions and different measures. For the latter, the representing
probability has density
\begin{equation}\label{eq:canonical-Wall-measure}
 \frac{d\sigma_{W_x}}{dt}(t)
 =\frac{\eta K_\eta(x,t)}
 {(1+t^2)|1-(x/\eta)F_x(t+i0)|^2}.
\end{equation}
To verify it, combine the formula for $\Im W_x$ with the Cauchy
boundary formula $\Im F_x(t+i0)=\pi\eta K_\eta(x,t)$ and
divide the Cauchy density by $1+t^2$.
The hypergeometric expression, or the differential equation away
from $\pm i$, continues the upper branch of $F_x$ analytically
through each real point. Its strictly positive boundary imaginary
part prevents the denominator in \eqref{eq:canonical-Wall-measure}
from vanishing there. Stieltjes inversion on each compact interval
therefore leaves no singular component on the real line. Since
$F_x(iy)\to0$, also $W_x(iy)\to0$, so the coefficient $a$ in
\eqref{eq:general-Nevanlinna} is zero. The normalization at $i$
then proves that \eqref{eq:canonical-Wall-measure} has mass one.

This identifies the measure of the canonical fraction as well as
that of the original Meixner fraction. It also explains why the
restriction $\eta>1$ belongs only to the unweighted probability
$\nu_x$: both normalized Nevanlinna probabilities exist for every
$\eta>0$.

\subsection{A joint law with Cauchy marginals}

There is a particularly simple way to put both spectral directions
in one probability space.

\begin{theorem}\label{thm:joint-expanded}
For every $\eta>0$, the density
\begin{equation}\label{eq:joint-expanded}
 \boxed{\displaystyle
 d\Pi_\eta(x,\delta)=
       \frac{K_\eta(x,\delta)}{\pi(1+\delta^2)}\,dx\,d\delta}
\end{equation}
defines a probability measure. Its marginals are
\begin{equation}\label{eq:joint-marginals}
 d\Pi_\eta^\delta(\delta)=\frac{d\delta}{\pi(1+\delta^2)},
 \qquad
 d\Pi_\eta^x(x)=\frac{\eta\,dx}{\pi(x^2+\eta^2)}.
\end{equation}
Its conditional laws are
\begin{equation}\label{eq:joint-conditionals}
 \Pi_\eta(dx\mid\delta)=d\mu_\delta(x),\qquad
 \Pi_\eta(d\delta\mid x)=d\widehat\sigma_x(\delta).
\end{equation}
\end{theorem}

\begin{proof}
Integrate first in $x$. The kernel has mass one, so the first
marginal and the total mass follow from Tonelli's theorem. Integrate
first in $\delta$ and use \eqref{eq:dual-weighted-mass} to get
the second marginal. Dividing the joint density by either marginal
gives the two conditional densities.
\end{proof}

This gives a literal answer to the question posed at the start.
The ordinary Meixner measure and the affine-normalized Nevanlinna
measure are the two conditional laws of one joint distribution.
The first marginal is standard Cauchy, and the second is Cauchy
with scale $\eta$. In the angle $\delta=\tan\theta$, the first
marginal is simply uniform on $(-\pi/2,\pi/2)$.

The unweighted probability $\nu_x$, when $\eta>1$, has a slightly
different role. It and $\mu_\delta$ are normalized slices of the
same kernel:
\[
 \frac{d\nu_x}{d\delta}=(\eta-1)K_\eta(x,\delta),\qquad
 \frac{d\mu_\delta}{dx}=K_\eta(x,\delta).
\]
The unweighted two-variable measure $K_\eta(x,\delta)\,dx\,d\delta$
has infinite total mass. The Cauchy reference factor in
\eqref{eq:joint-expanded} turns it into a probability and converts
the second slice into $\widehat\sigma_x$.
Orthogonality in each variable does not by itself supply a joint
probability law; the Cauchy reference factor supplies the required
compatibility.

As a useful check, at $\eta=2$ and $x=0$,
\[
 K_2(0,t)=\frac1{\pi(1+t^2)},\qquad
 F_0(z)=-\frac2{z+i},\qquad
 d\widehat\sigma_0(t)=\frac{2\,dt}{\pi(1+t^2)^2}.
\]
The Cauchy measure $2K_2(0,t)\,dt$ has mass two, while the finite
Nevanlinna measure has mass one. This small example makes an
omitted factor $1+t^2$ impossible to overlook. In the other
direction, the same kernel gives the tangent density
\eqref{eq:tangent-density}.

\subsection{The parameter law as a conjugate posterior}

The angle also gives the joint law a familiar statistical meaning.
For a fixed parameter $\theta$, the density of the observation $x$ is
\[
 K_\eta(x,\tan\theta)
       =C_\eta(x)e^{x\theta}\cos^\eta\theta.
\]
The dependence on $\theta$ is entirely through an exponential and
a power of the cosine. It is therefore natural to give $\theta$
a prior density of the same form,
\begin{equation}\label{eq:conjugate-prior}
 p_{r,s}(\theta)=\frac{e^{s\theta}\cos^r\theta}{Z(r,s)},
 \qquad -\frac\pi2<\theta<\frac\pi2,\quad r>-1,\quad s\in\R.
\end{equation}
The normalizing constant is
\begin{equation}\label{eq:Z}
 Z(r,s)=\int_{-\pi/2}^{\pi/2}e^{s\theta}\cos^r\theta\,d\theta
 =\frac{\pi\Gamma(r+1)}
 {2^r\Gamma(1+(r+is)/2)\Gamma(1+(r-is)/2)}.
\end{equation}
No new integral is needed: substitute $t=\tan\theta$ into
\eqref{eq:dual-beta} and take $a=1+(r+is)/2$. Its condition
$\Re a>1/2$ becomes exactly $r>-1$.

Multiplication by the likelihood merely adds exponents. Bayes'
formula gives the posterior update and the predictive density
\begin{equation}\label{eq:bayes-update}
 (r,s)\longmapsto(r+\eta,s+x),\qquad
 p_{\eta;r,s}(x)
       =C_\eta(x)\frac{Z(r+\eta,s+x)}{Z(r,s)}.
\end{equation}
Thus the family of priors is \emph{conjugate}: observing $x$ keeps
the posterior within the same two-parameter family. For $r=s=0$
the prior is uniform in angle and hence standard Cauchy in
$\delta$. Formula \eqref{eq:dual-weighted-mass}, or gamma
cancellation in \eqref{eq:Z}, gives the scale-$\eta$ Cauchy
predictive density. The posterior, expressed in $\delta$, is
exactly $\widehat\sigma_x$ in \eqref{eq:phi-probability}.
This conjugate family and its predictive distributions are discussed
in \cite{DevroyeHill2025,TojoYoshino2026}. Here they explain why
the two spectral measures fit so simply into one joint probability.

\subsection{Bounded rational moments}

Under $\widehat\sigma_x$, the angle has density
$e^{x\theta}\cos^\eta\theta/Z(\eta,x)$. Its Fourier coefficients
have a particularly simple product form. For every integer $n\ge0$,
\begin{equation}\label{eq:cayley-moments}
 \boxed{\displaystyle
 \int_{\R}\left(\frac{1+it}{1-it}\right)^n
          d\widehat\sigma_x(t)
 =\prod_{j=0}^{n-1}
       \frac{\eta-2j+ix}{\eta+2j+2-ix}.}
\end{equation}
The empty product is one, and negative powers give the complex
conjugates. To prove the formula, observe that
\[
 e^{2i\theta}=\frac{1+it}{1-it},\qquad t=\tan\theta.
\]
The required integral is therefore $Z(\eta,x+2in)/Z(\eta,x)$.
For a fixed $r>-1$, both sides of \eqref{eq:Z} are entire in $s$:
on the integral side the interval is bounded, and on the gamma
side one uses the entire reciprocal gamma function. The identity
for real $s$ therefore extends to complex $s$. Gamma recurrence
now gives the displayed product; if a reciprocal gamma factor
vanishes, the result follows by continuity.

The integrand in \eqref{eq:cayley-moments} has modulus one on
the real line. These rational moments exist at every degree even
though sufficiently high polynomial moments do not. The circle
coordinate is precisely the Cayley coordinate of
\eqref{eq:Wall-Cayley}: what looked like a difficulty at infinity
on the line becomes ordinary Fourier analysis on the circle.
This is the natural setting of the complementary
Romanovski--Routh and unit-circle connections
\cite{MartinezFinkelshteinEtAl2019,IsmailSriRanga2019}.


\section{A joint Cauchy law at every finite degree}\label{sec:finite}

The joint law of the preceding section is already present in every
finite continued fraction. At finite degree it lives on the zero
curves of a polynomial, yet its two marginals are exactly the same
Cauchy laws as in the limit. To see why, we need only follow a zero
as the parameter moves.

\subsection{The zeros form increasing spectral curves}

Return to the positive pencil of
\eqref{eq:two-variable-determinant}, and write
\[
 H_N(\delta)=B_N+\delta A_N,\qquad
 x_1(\delta)<\cdots<x_N(\delta)
\]
for its eigenvalues. They are simple because a diagonal unitary
conjugation makes $H_N(\delta)$ an irreducible real Jacobi matrix.
If $v_j(\delta)$ is a unit eigenvector, differentiating its eigenvalue
equation and taking the inner product with $v_j$ gives the
Hellmann--Feynman formula
\begin{equation}\label{eq:velocity}
 x_j'(\delta)=v_j(\delta)^*A_Nv_j(\delta)>0.
\end{equation}
Positivity of the pencil has become motion in a definite direction.

There are exact, familiar bounds for this motion. The characteristic
polynomials of the leading blocks of $A_N$ satisfy
\[
 \mathcal L_{n+1}(u)
 =(u-2n-\eta)\mathcal L_n(u)
       -n(n+\eta-1)\mathcal L_{n-1}(u),
 \qquad \mathcal L_0=1,\quad\mathcal L_1=u-\eta.
\]
These are the monic Laguerre polynomials
$\mathcal L_n(u)=(-1)^nn!L_n^{(\eta-1)}(u)$
\cite[\S9.12]{KoekoekLeskySwarttouw2010}.
Thus, if $0<\ell_1<\cdots<\ell_N$ are the zeros of
$L_N^{(\eta-1)}$, the Rayleigh quotient in \eqref{eq:velocity} gives
\begin{equation}\label{eq:velocity-bounds}
 \ell_1\le x_j'(\delta)\le\ell_N.
\end{equation}
Each zero curve is consequently a real analytic increasing bijection
of $\R$ onto $\R$. The two spectral problems are the vertical and
horizontal readings of the same $N$ disjoint curves.

The pencil gives a two-variable statement as well. A real polynomial
is called \emph{real stable} if it is nonzero whenever all its
variables have positive imaginary parts. Here
$p_N(x,y)=M_N(x;-y,\eta)$ is real stable: the imaginary part of its
matrix $xI_N+yA_N-B_N$ is
$(\Im x)I_N+(\Im y)A_N>0$, so that matrix cannot have a null vector.
More generally, every oblique slice
\begin{equation}\label{eq:oblique}
 t\longmapsto M_N(a+ut;b-vt,\eta),
 \qquad a,b\in\R,\quad u,v\ge0,\quad u+v>0,
\end{equation}
has only real zeros. Its matrix is a Hermitian pencil with positive
coefficient $uI_N+vA_N$; conjugating by the inverse square root of
that coefficient turns its zeros into ordinary Hermitian eigenvalues.
The reality of the zeros in either variable is part of this more
coherent geometric fact.

\subsection{A Jacobian relates the two quadrature weights}

The spectral weights in the $x$ direction are
\begin{equation}\label{eq:ordinary-quadrature}
 w_j(\delta)=|\langle e_0,v_j(\delta)\rangle|^2,\qquad
 \mu_{\delta,N}=\sum_{j=1}^Nw_j(\delta)\,
                         \delta_{x_j(\delta)}.
\end{equation}
Here $\delta_a$ denotes a unit point mass at $a$.
This is the Gaussian quadrature measure used in
Section~\ref{sec:analytic-fractions}: the diagonal unitary change
of basis preserves the first-coordinate weights. Each weight is
positive, and their sum is one.

For the parameter direction, factor the matrix in
\eqref{eq:Wall-finite-resolvent} through $A_N^{1/2}$. If
\[
 T_{x,N}=A_N^{-1/2}(xI_N-B_N)A_N^{-1/2},
\]
the spectral measure of $T_{x,N}$ at the vector
$\sqrt\eta\,A_N^{-1/2}e_0$ is a positive discrete measure
$\rho_{x,N}$, and
\begin{equation}\label{eq:finite-F}
 F_{x,N}(z)=\int_\R\frac{\rho_{x,N}(\dd t)}{t-z}.
\end{equation}
Its support consists of the parameter zeros of $M_N(x;t,\eta)$.

\begin{proposition}\label{prop:jacobian}
At the crossing $x=x_j(t)$, the two weights satisfy
\begin{equation}\label{eq:jacobian}
 \rho_{x,N}(\{t\})=\frac{\eta w_j(t)}{x_j'(t)}.
\end{equation}
Consequently, as locally finite measures on $\R^2$,
\begin{equation}\label{eq:coarea}
 \mu_{\delta,N}(\dd x)\dd\delta
 =\frac1\eta\rho_{x,N}(\dd\delta)\dd x.
\end{equation}
\end{proposition}

\begin{proof}
At a crossing, the corresponding unit eigenvector of $T_{x,N}$ is
\[
 q=\frac{A_N^{1/2}v_j(t)}
          {\sqrt{v_j(t)^*A_Nv_j(t)}}.
\]
Its spectral weight at $\sqrt\eta\,A_N^{-1/2}e_0$ is
$\eta w_j(t)/(v_j(t)^*A_Nv_j(t))$. Formula
\eqref{eq:velocity} proves \eqref{eq:jacobian}.

Now substitute $x=x_j(\delta)$ on each increasing curve.
For a compactly supported continuous function $f$,
\[
 \sum_j\int_\R f(x_j(\delta),\delta)w_j(\delta)\dd\delta
 =\int_\R\sum_{j=1}^N
  f(x,t_j(x))\frac{w_j(t_j(x))}{x_j'(t_j(x))}\dd x,
\]
where $t_j=x_j^{-1}$. This is \eqref{eq:coarea}.
\end{proof}

The Jacobian is the essential point. The two finite measures assign
mass along the same curve, using different coordinates on that curve.

\subsection{The finite coupling and its limit}

\begin{theorem}\label{thm:finite-coupling}
For every $N\ge1$ and $\eta>0$, the probability measure
\begin{equation}\label{eq:finite-joint}
 \Pi_{\eta,N}(\dd x,\dd\delta)
 =\mu_{\delta,N}(\dd x)\frac{\dd\delta}{\pi(1+\delta^2)}
\end{equation}
is supported on $M_N(x;\delta,\eta)=0$ and has the exact marginals
\begin{equation}\label{eq:finite-marginals}
 \Pi_{\eta,N}^{\delta}(\dd\delta)
   =\frac{\dd\delta}{\pi(1+\delta^2)},\qquad
 \Pi_{\eta,N}^{x}(\dd x)
   =\frac{\eta\dd x}{\pi(x^2+\eta^2)}.
\end{equation}
Its reverse conditional law is
\begin{equation}\label{eq:finite-reverse}
 \widehat\sigma_{x,N}(\dd t)
 =\frac{x^2+\eta^2}{\eta^2(1+t^2)}\rho_{x,N}(\dd t).
\end{equation}
As $N\to\infty$, $\Pi_{\eta,N}$ converges weakly to the joint
law $\Pi_\eta$ of \eqref{eq:joint-expanded}.
\end{theorem}

\begin{proof}
The first marginal follows from $\mu_{\delta,N}(\R)=1$.
The common value $F_{x,N}(i)=F_x(i)=\eta/(x-i\eta)$, obtained
in the proof of Theorem~\ref{thm:Wall-fraction-expanded}, gives
\[
 \int_\R\frac{\rho_{x,N}(\dd t)}{1+t^2}
 =\Im F_{x,N}(i)=\frac{\eta^2}{x^2+\eta^2}.
\]
Using \eqref{eq:coarea}, the $x$-marginal density of
\eqref{eq:finite-joint} is therefore
\[
 \frac1{\pi\eta}\int_\R\frac{\rho_{x,N}(\dd t)}{1+t^2}
 =\frac{\eta}{\pi(x^2+\eta^2)}.
\]
Dividing the same measure identity by this density proves
\eqref{eq:finite-reverse}.

For fixed $\delta$, the convergence
$\mu_{\delta,N}\Rightarrow\mu_\delta$ was proved in
Theorem~\ref{thm:Jacobi-fraction}. Apply it to
$f(\,\cdot\,,\delta)$, for a bounded continuous function
$f$ on $\R^2$, and then use dominated convergence against the fixed
Cauchy law in $\delta$. This proves the joint weak convergence.
\end{proof}

\begin{figure}[!htbp]
 \centering
 \includegraphics[width=\linewidth]{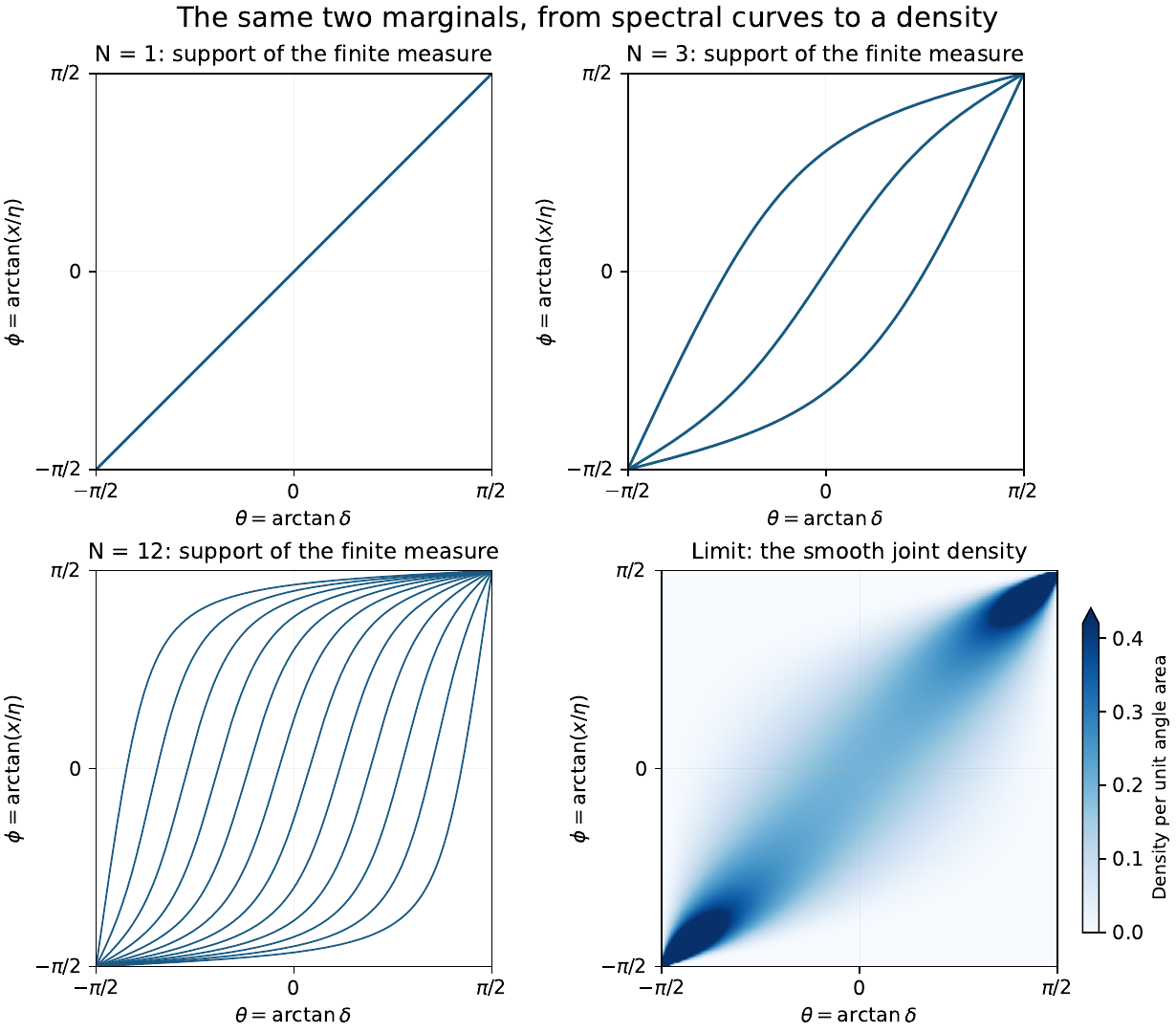}
 \caption{Finite spectral curves and the limiting density, with $\eta=2$,
 in the angles $\theta=\arctan\delta$ and $\varphi=\arctan(x/\eta)$.
 Both angular marginals are uniform on $(-\pi/2,\pi/2)$ at every stage.
 The first three panels show support only: the conditional weights are
 those in \eqref{eq:ordinary-quadrature}, and are not encoded by line
 thickness. The last panel shows the joint density in the two angles.}
 \label{fig:spectral-curves}
\end{figure}

For $N=1$, the coupling is simply $x=\eta\delta$.
Higher degrees change the dependence while preserving both marginals.
The smooth density of \eqref{eq:joint-expanded} is approached by
singular measures on algebraic curves, as Figure~\ref{fig:spectral-curves}
illustrates.

\clearpage

This is spectral averaging in a particularly explicit form; general
matrix-valued spectral averaging is developed in
\cite{LiawTreil2020}. Here the positive pencil makes the mechanism
visible: a zero moves with velocity $v^*A_Nv$, and its spectral mass
changes by the reciprocal velocity when the two coordinates are
flipped. The finite fractions therefore retain both measures
and the compatibility between them.


\section{What the finite fractions tell us quantitatively}\label{sec:quantitative-fractions}

The finite measures also give practical formulas and convergence
bounds. The general relationship between $R_{\mathrm{II}}$
recurrences and rational quadrature is established
\cite{BraccialiEtAl2020}; here all the constants are explicit.

\subsection{Positive rational Gaussian weights}

The preceding section identified the finite measures geometrically.
We can now write each of their masses in closed form. Recall that
$F_{x,N}=\eta N_N/M_N$, with numerator polynomials $N_0=0$, $N_1=1$
from the continuant recurrence in Section~\ref{sec:parameter-duality}.
The exact identity \eqref{eq:Wall-exact-error} gives
\[
 F_x(z)-F_{x,N}(z)=O((z-i)^N)\qquad(z\to i).
\]
This interpolation property will provide both quadrature exactness
and a convergence bound.

\begin{proposition}\label{prop:quadrature}
Let $t_1<\cdots<t_N$ be the parameter zeros of
$M_N(x;t,\eta)$, and let a prime denote differentiation in $t$.
Then
\begin{align}
 \rho_{x,N}(\{t_j\})
 &=\frac{\eta N!(\eta)_N(1+t_j^2)^{N-1}}
        {[M_N'(x;t_j,\eta)]^2},\label{eq:rho-weights}\\
 \omega_j:=\widehat\sigma_{x,N}(\{t_j\})
 &=\frac{x^2+\eta^2}{\eta}
   \frac{N!(\eta)_N(1+t_j^2)^{N-2}}
        {[M_N'(x;t_j,\eta)]^2}.\label{eq:prob-weights}
\end{align}
All weights are positive and $\sum_j\omega_j=1$. For every
polynomial $p$ of degree at most $2N-1$,
\begin{equation}\label{eq:rational-quadrature}
 \int_\R\frac{p(t)}{(1+t^2)^{N-1}}\widehat\sigma_x(\dd t)
 =\sum_{j=1}^N\omega_j\frac{p(t_j)}{(1+t_j^2)^{N-1}}.
\end{equation}
These statements hold for every $\eta>0$.
\end{proposition}

\begin{proof}
The recurrence gives the continuant identity
\[
 M_{N-1}N_N-M_NN_{N-1}
 =(N-1)!(\eta)_{N-1}(1+t^2)^{N-1}.
\]
The parameter lowering identity \eqref{eq:delta-lowering-expanded}
then gives
\[
 N_N(x;t_j,\eta)
 =-\frac{N!(\eta)_N(1+t_j^2)^{N-1}}{M_N'(x;t_j,\eta)}.
\]
Taking the negative residue of $F_{x,N}=\eta N_N/M_N$
proves \eqref{eq:rho-weights}; \eqref{eq:finite-reverse}
proves \eqref{eq:prob-weights}.

Taylor interpolation at $i$ gives equality of the integrals
of $(t-i)^{-k}$ against $\rho_x$ and $\rho_{x,N}$ for
$1\le k\le N$. Reflection gives the same assertion at $-i$.
Every $p(t)/(1+t^2)^N$ with $\deg p\le2N-1$ is a linear
combination of those $2N$ partial fractions. Divide the measures
by $1+t^2$ and use their common normalization to obtain
\eqref{eq:rational-quadrature}. All integrals converge: at the
largest allowed degree the integrand against $\rho_x$ is
$O(|t|^{-\eta-1})$.
\end{proof}

This is genuine Gaussian exactness in a rational space. The
denominator provides the moments which the heavy-tailed measure
does not possess as ordinary polynomial moments.

\subsection{The mass transition at $\eta=1$}

The total mass of a finite Wall measure is independent of $x$:
\begin{equation}\label{eq:finite-mass}
 \rho_{x,N}(\R)=\eta(A_N^{-1})_{00}
 =\eta\sum_{k=0}^{N-1}\frac{k!}{(\eta)_{k+1}}
 =\begin{cases}
 \displaystyle\frac{\eta}{\eta-1}
       \left(1-\frac{N!}{(\eta)_N}\right),&\eta\ne1,\\[6pt]
 H_N,&\eta=1.
 \end{cases}
\end{equation}
Here $H_N=\sum_{j=1}^N1/j$ is the $N$th harmonic number.
To prove the formula, solving $L_Nu=e_0$ gives
$|u_k|^2=k!/(\eta)_{k+1}$, and the sum telescopes because
\[
 \frac{k!}{(\eta)_k}-\frac{(k+1)!}{(\eta)_{k+1}}
       =(\eta-1)\frac{k!}{(\eta)_{k+1}}.
\]
Thus the finite measures approach total mass $\eta/(\eta-1)$
when $\eta>1$, grow as $\log N$ when $\eta=1$, and grow as
\[
 \frac{\eta\Gamma(\eta)}{1-\eta}N^{1-\eta}
 \qquad(0<\eta<1).
\]
The threshold for the infinite measure is already encoded in the
first column of a finite triangular inverse. The normalized
probabilities \eqref{eq:finite-reverse} have mass one throughout.

\subsection{A geometric error bound for the Wall fraction}

\begin{proposition}\label{prop:error}
For $z\in\C_+$, set
\[
 q=\left|\frac{z-i}{z+i}\right|<1,
 \qquad \alpha=\frac{\eta^2}{x^2+\eta^2}.
\]
Then
\begin{equation}\label{eq:error}
 |F_{x,N}(z)-F_x(z)|
 \le\frac{4\alpha q^N}{(1-q)^2}.
\end{equation}
\end{proposition}

\begin{proof}
Put $c=F_x(i)=F_{x,N}(i)$, so $\Im c=\alpha$, and map the
upper half-plane to the disk by $w=(z-i)/(z+i)$. The functions
\[
 S(w)=\frac{F_x(z)-c}{F_x(z)-\overline c},\qquad
 S_N(w)=\frac{F_{x,N}(z)-c}{F_{x,N}(z)-\overline c}
\]
map the disk into itself and vanish at zero. Schwarz's lemma
gives $|S|,|S_N|\le|w|$. Their difference vanishes to order
$N$ by \eqref{eq:Wall-exact-error}; the same lemma applied to
$(S_N-S)/2$ gives $|S_N-S|\le2|w|^N$.
Inverting the fractional linear transformation yields
\[
 |F_{x,N}-F_x|
 =\frac{2\alpha|S_N-S|}{|1-S_N||1-S|},
\]
which proves the bound.
\end{proof}

The interpolation point $i$ is the natural center of this estimate:
its Cayley coordinate $q$ measures how far we have moved toward the
boundary. On every compact subset of the upper half-plane the
continued fraction therefore converges geometrically. No asymptotic
analysis of the hypergeometric expression is required.


\section{What the two Cauchy variables do at infinity}
\label{sec:research-tails}

The two marginals remain Cauchy while the shape parameter $\eta$
changes. What changes is the dependence between the coordinates.
Ordinary correlation cannot measure it, since the Cauchy marginals
have no second moments. We can ask a more direct question: if one
coordinate is large and positive, how likely is the other to be
comparably large? The answer leads back to the gamma measure and
the Laguerre polynomials that bounded the velocities of our zeros.

Rescale the first coordinate by putting $Y=X/\eta$ and write $D$ for
the parameter coordinate. Both $D$ and $Y$ now have the standard Cauchy
law
\[
 m(\dd d)=\frac{\dd d}{\pi(1+d^2)}.
\]
Their conditional kernel and joint density are
\begin{equation}\label{eq:research-Q-kernel}
 Q_\eta(d,\dd y)=\eta K_\eta(\eta y,d)\,\dd y,
 \qquad
 p_\eta(d,y)=\frac{\eta K_\eta(\eta y,d)}{\pi(1+d^2)}.
\end{equation}
The integral operator
$Q_\eta f(d)=\int f(y)Q_\eta(d,\dd y)$ is the $1+d^2$ member of the
Ismail--May family of approximation operators
\cite{IsmailMay1978}; its kernel appears explicitly in
\cite[Eq.~(1)]{GuptaMilovanovic2022}. Our concern here is its behavior
under its invariant Cauchy measure.

The Meixner transform gives
\begin{equation}\label{eq:research-Q-transform}
 \E[e^{itY}\mid D=d]
 =\bigl[\cosh(t/\eta)-id\sinh(t/\eta)\bigr]^{-\eta},
 \qquad
 \E[Y\mid D=d]=d,\quad
 \operatorname{Var}(Y\mid D=d)=\frac{1+d^2}{\eta}.
\end{equation}
The complex power is the continuous branch equal to one at $t=0$.
The conditional moments here are moments of the probability
$Q_\eta(d,\dd y)$ at a fixed $d$. They are finite even though the
stationary pair $(D,Y)$ has no integrable coordinates.

\subsection{A gamma law at the edge}

\begin{proposition}\label{prop:research-gamma-ratio}
Fix $\eta>0$. As $d\to+\infty$,
\begin{equation}\label{eq:research-gamma-ratio}
 \mathcal L(Y/d\mid D=d)\Longrightarrow G_\eta,
 \qquad
 \Pbb(G_\eta\in\dd u)
 =\frac{\eta^\eta}{\Gamma(\eta)}u^{\eta-1}e^{-\eta u}\,\dd u,
 \quad u>0.
\end{equation}
The same ratio limit holds as $d\to-\infty$.
\end{proposition}

\begin{proof}
Replace $t$ by $t/d$ in \eqref{eq:research-Q-transform}. Then
\[
 \E[e^{itY/d}\mid D=d]
 =\left[\cosh\!\left(\frac{t}{\eta d}\right)
   -id\sinh\!\left(\frac{t}{\eta d}\right)\right]^{-\eta}
 \longrightarrow(1-it/\eta)^{-\eta}.
\]
The limit is the characteristic function of the displayed gamma law.
L\'evy's continuity theorem proves the assertion. Reflection,
$K_\eta(-x,-d)=K_\eta(x,d)$, gives the other endpoint.
\end{proof}

The limiting relative fluctuation has mean one and variance
$1/\eta$. Its gamma law is the Laguerre orthogonality measure, with
the scale changed to give mean one.

Reversing the conditioning gives another classical distribution.
It is useful to state the reverse limit separately, since it explains
why the two directions have different normalizations.

\begin{proposition}\label{prop:inverse-gamma-ratio}
As $y\to+\infty$, the conditional law of $D/y$ given $Y=y$
converges in total variation to the inverse-gamma law with density
\begin{equation}\label{eq:inverse-gamma-ratio}
 \frac{\eta^{\eta+1}}{\Gamma(\eta+1)}
 v^{-\eta-2}e^{-\eta/v},\qquad v>0.
\end{equation}
The same assertion holds as $y\to-\infty$.
\end{proposition}

\begin{proof}
The reverse conditional law \eqref{eq:shifted-eta-measure} is
\[
 \Pbb(D\in\dd d\mid Y=y)
 =(\eta+1)K_{\eta+2}(\eta y,d)\,\dd d.
\]
For positive $y$, the density of $D/y$ is therefore
$y(\eta+1)K_{\eta+2}(\eta y,yv)$. Use the gamma asymptotic
\eqref{eq:Ceta-tail} and, for $v>0$,
\[
 \arctan(yv)=\frac\pi2-\frac1{yv}+O(y^{-3}).
\]
Substitution gives \eqref{eq:inverse-gamma-ratio} pointwise.
For $v<0$ the two exponential factors instead leave a factor
$e^{-\pi\eta y+O(1)}$, so the limit is zero.
The limiting density has mass one, as the substitution
$u=1/v$ shows. Scheff\'e's lemma proves total-variation
convergence. Central reflection gives the other endpoint.
\end{proof}

This inverse-gamma variable is the reciprocal of a
\emph{size-biased} $G_\eta$: since $\E G_\eta=1$, size bias means
multiplying its density by $u$, which changes its shape from
$\eta$ to $\eta+1$. Thus changing the direction of conditioning
at the boundary has a precise and elementary cost.

\subsection{An exact answer for simultaneous extremes}

\begin{theorem}\label{thm:research-joint-tail}
For every $\eta>0$ and $a,b>0$,
\begin{equation}\label{eq:research-joint-tail}
 \boxed{\displaystyle
 \lim_{R\to\infty}\pi R\,\Pbb(D>aR,\ Y>bR)
 =\E\min\!\left(\frac1a,\frac{G_\eta}{b}\right).}
\end{equation}
The same formula holds for the two lower tails. Opposite-sign tails
satisfy
\begin{equation}\label{eq:research-opposite-tail}
 \lim_{R\to\infty}\pi R\,\Pbb(D>aR,\ Y<-bR)=0.
\end{equation}
\end{theorem}

\begin{proof}
Disintegrate over $D$ and put $d=Rr$:
\[
 \pi R\,\Pbb(D>aR,Y>bR)
 =\int_a^\infty\frac{R^2}{1+R^2r^2}
  \Pbb\!\left(\frac{Y}{Rr}>\frac br\,\middle|\,D=Rr\right)\dd r.
\]
Proposition~\ref{prop:research-gamma-ratio} identifies the pointwise
limit. The integrand is bounded by $r^{-2}$, which is integrable on
$[a,\infty)$. Dominated convergence and Tonelli's theorem give
\begin{align*}
 \lim_{R\to\infty}\pi R\,\Pbb(D>aR,Y>bR)
 &=\int_a^\infty r^{-2}\Pbb(G_\eta>b/r)\,\dd r\\
 &=\E\left[\int_a^\infty r^{-2}
       \mathbf1_{\{r>b/G_\eta\}}\,\dd r\right]
 =\E\min(1/a,G_\eta/b).
\end{align*}
For opposite signs the limiting probability is zero, since
$G_\eta>0$. Reflection proves the lower-tail assertion.
\end{proof}

For an explicit unequal-threshold formula, let
$\mathsf P(s,z)=\gamma(s,z)/\Gamma(s)$ and
$\mathsf Q(s,z)=1-\mathsf P(s,z)$ denote the regularized incomplete
gamma functions. Splitting the expectation at $G_\eta=b/a$ gives
\begin{equation}\label{eq:research-tail-incomplete-gamma}
 \E\min(1/a,G_\eta/b)
 =\frac1a\mathsf Q\!\left(\eta,\frac{\eta b}{a}\right)
  +\frac1b\mathsf P\!\left(\eta+1,\frac{\eta b}{a}\right).
\end{equation}

\begin{corollary}\label{cor:research-tail-coefficient}
The probability that $Y$ exceeds the same high threshold as $D$ has
the limit
\begin{equation}\label{eq:research-tail-coefficient}
 \boxed{\displaystyle
 \lambda_\eta:=\lim_{R\to\infty}\Pbb(Y>R\mid D>R)
 =1-\frac{\eta^\eta e^{-\eta}}{\Gamma(\eta+1)}.}
\end{equation}
It is strictly increasing in $\eta$, from zero to one.
\end{corollary}

\begin{proof}
Use $\pi R\Pbb(D>R)\to1$ and set $a=b=1$ in
\eqref{eq:research-tail-incomplete-gamma}. The incomplete-gamma
recurrence reduces the answer to \eqref{eq:research-tail-coefficient}.

For monotonicity, take $0<\alpha<\beta$ and independent gamma variables
$A\sim\Gamma(\alpha,1)$, $B\sim\Gamma(\beta-\alpha,1)$, where the
second parameter is the rate. Changing variables from $(A,B)$
to $(A+B,A/(A+B))$ in their product density shows that the ratio
is independent of the sum and has beta mean $\alpha/\beta$.
Consequently,
\[
 \E[A/\alpha\mid A+B]=(A+B)/\beta.
\]
Apply conditional Jensen to the concave function $u\mapsto\min(1,u)$
and take expectations:
$\E\min(1,A/\alpha)\le\E\min(1,(A+B)/\beta)$.
The inequality is strict: when $A+B>\alpha$, the conditional law of
$A/\alpha$ places positive mass on both sides of its corner at one.
The endpoint limits follow from the displayed formula.
\end{proof}

For example, $\lambda_1=1-e^{-1}$. More precisely, the gamma expansion
at zero and Stirling's formula give
\begin{align}
 \lambda_\eta
 &=\eta\log(1/\eta)+(1-\gamma)\eta
   +O\!\left(\eta^2\log^2\eta\right)
 &&(\eta\downarrow0),\label{eq:research-tail-small-eta}\\
 1-\lambda_\eta
 &=\frac1{\sqrt{2\pi\eta}}
   \left(1-\frac1{12\eta}+O(\eta^{-2})\right)
 &&(\eta\to\infty).\label{eq:research-tail-large-eta}
\end{align}
Here $\gamma$ is Euler's constant. These finite quantities describe
dependence even though a covariance cannot do so.

\subsection{From independence to coincidence}

\begin{proposition}\label{prop:research-shape-endpoints}
As $\eta\downarrow0$, the law of $(D,Y)$ tends in total variation to
the product of two standard Cauchy laws. As $\eta\to\infty$, it tends
weakly to the diagonal law $(D,D)$, with $D\sim m$.
\end{proposition}

\begin{proof}
At the first endpoint, $\Gamma(z)\sim1/z$ gives, for fixed $y$,
\[
 \eta C_\eta(\eta y)\longrightarrow\frac1{\pi(1+y^2)}.
\]
The remaining factor
$(1+d^2)^{-\eta/2}e^{\eta y\arctan d}$ tends to one. Hence
$p_\eta(d,y)$ converges pointwise to the product Cauchy density;
Scheff\'e's lemma gives total-variation convergence.

At the other endpoint, \eqref{eq:research-Q-transform} and Chebyshev's
inequality give, for each $\varepsilon>0$,
\[
 \Pbb(|Y-d|>\varepsilon\mid D=d)
 \le\min\!\left(1,\frac{1+d^2}{\eta\varepsilon^2}\right).
\]
Integrating against $m$ and applying dominated convergence shows that
$Y-D\to0$ in probability. This proves the weak limit using only the
bounded conditional probability estimate.
\end{proof}

\subsection{Laguerre quadrature at the boundary}

The finite couplings on spectral curves have the same boundary mechanism.
Let $\ell_{j,N}$ and $\alpha_{j,N}$ be the nodes and normalized weights
of $N$-point Gaussian quadrature for
$u^{\eta-1}e^{-u}\,\dd u/\Gamma(\eta)$. Equivalently, they are the
eigenvalues and first-coordinate spectral weights of $A_N$.

\begin{corollary}\label{cor:research-finite-tail}
Under $\Pi_{\eta,N}$ put $Y_N=X/\eta$, and let
\[
 G_{\eta,N}\sim\sum_{j=1}^N
       \alpha_{j,N}\delta_{\ell_{j,N}/\eta}.
\]
For $a,b>0$,
\begin{equation}\label{eq:research-finite-joint-tail}
 \lim_{R\to\infty}\pi R\,\Pbb(D>aR,Y_N>bR)
 =\E\min(1/a,G_{\eta,N}/b).
\end{equation}
In particular,
\begin{equation}\label{eq:research-finite-tail-coefficient}
 \begin{aligned}
 \lambda_{\eta,N}
 &:=\lim_{R\to\infty}\Pbb(Y_N>R\mid D>R)\\
 &=\sum_{j=1}^N\alpha_{j,N}\min(1,\ell_{j,N}/\eta)
 \longrightarrow\lambda_\eta\qquad(N\to\infty).
 \end{aligned}
\end{equation}
\end{corollary}

\begin{proof}
As $d\to+\infty$, we have
$d^{-1}H_N(d)=A_N+d^{-1}B_N\to A_N$ in matrix norm. Continuity
of simple eigenvalues and their spectral projections gives
\[
 x_j(d)/d\to\ell_{j,N},\qquad
 w_j(d)\to\alpha_{j,N}.
\]
Thus $Y_N/d$, conditional on $D=d$, converges to $G_{\eta,N}$.
Repeat the proof of Theorem~\ref{thm:research-joint-tail}. The limiting
conditional probability can jump only when $b/r$ equals one of its
finitely many atoms; these exceptional $r$'s do not affect dominated
convergence.

The laws of $G_{\eta,N}$ are the Gaussian quadrature measures for
$G_\eta$ and converge weakly to it by the same moment argument used
in Theorem~\ref{thm:Jacobi-fraction}: their common mean gives
tightness on $[0,\infty)$, higher matched moments give uniform
integrability for each fixed power, and the gamma law is determined
by its moment generating function near zero. Finally,
$u\mapsto\min(1,u)$ is bounded and continuous.
\end{proof}

The first quadrature has
$G_{\eta,1}=1$, hence $\lambda_{\eta,1}=1$. The convergence in
\eqref{eq:research-finite-tail-coefficient} follows from the boundary
quadrature calculation; weak convergence of the joint measures alone
would not imply convergence of their tail coefficients.

For a small explicit check, the two Laguerre nodes are
$\eta+1\mp\sqrt{\eta+1}$, with weights
$(1\pm1/\sqrt{\eta+1})/2$. They give
\[
 \lambda_{\eta,2}=1-\frac1{2\sqrt{\eta+1}}.
\]
There is a second useful connection with the finite fraction:
\begin{equation}\label{eq:tail-inverse-mass}
 \E\frac1{G_{\eta,N}}
 =\eta\sum_j\frac{\alpha_{j,N}}{\ell_{j,N}}
 =\eta(A_N^{-1})_{00}=\rho_{x,N}(\R).
\end{equation}
The finite Wall mass is exactly the reciprocal moment of its gamma
quadrature. In the limit, the gamma reciprocal moment equals
$\eta/(\eta-1)$ for $\eta>1$ and is infinite for $\eta\le1$.
This gives the threshold in \eqref{eq:finite-mass} a probabilistic
meaning. The explicit sum there, rather than weak convergence alone,
justifies the convergence or divergence of these unbounded moments.

\subsection{Why the conditional transform is not compact}

The boundary gamma law also answers an operator-theoretic question.
Invariance of $m$ and conditional Jensen imply
\[
 \|Q_\eta f\|_{L^2(m)}^2
 \le\int Q_\eta|f|^2\,\dd m=\|f\|_{L^2(m)}^2.
\]
Thus $Q_\eta$ preserves nonnegative functions and is a contraction.
Its smooth kernel might
suggest compactness, but its behavior far out in the tails prevents it.

\begin{proposition}\label{prop:research-noncompact}
For every $\eta>0$, $Q_\eta$ is not compact on $L^2(m)$.
\end{proposition}

\begin{proof}
Set
\[
 f_R=\frac{\mathbf1_{[R,2R]}}{\sqrt{m([R,2R])}}.
\]
These functions have norm one and converge weakly to zero: for
$h\in L^2(m)$,
\[
 |\langle f_R,h\rangle|
 \le\left(\int_R^\infty|h|^2\,\dd m\right)^{1/2}\longrightarrow0.
\]
Put $c_\eta=\frac12\Pbb(1<G_\eta<4/3)>0$.
Proposition~\ref{prop:research-gamma-ratio} shows that, for all
sufficiently large $d$,
\[
 \Pbb(1<Y/d<4/3\mid D=d)\ge c_\eta.
\]
If $R\le d\le3R/2$, this event implies $R<Y<2R$. Therefore
\[
 \|Q_\eta f_R\|_2^2
 \ge c_\eta^2\frac{m([R,3R/2])}{m([R,2R])}
 \longrightarrow\frac23c_\eta^2>0.
\]
A compact operator would send this bounded weakly null sequence to a
norm-null sequence, a contradiction.
\end{proof}

In particular, the transform is not Hilbert--Schmidt. The two
conditional laws therefore provide a useful positive transform and
its adjoint without automatically supplying a compact spectral
decomposition.

The two measures, their finite quadratures, and their behavior at
infinity have now emerged from the same positive kernel. We can
return to the word \emph{bispectral} with a clearer picture of what
each spectral direction means.


\section{What is spectral, and what is bispectral?}\label{sec:spectral-meaning}

The two-variable kernel is concrete enough that we can distinguish
several structures often gathered under the word ``bispectrality.''
There is the classical relation between degree and argument, the
parameter duality of the positive pencil, and a representation-theoretic
explanation for the recurrence itself. Each answers a different question.

\subsection{The classical difference operator}

In the argument $x$, define the fixed difference operator
\begin{align}\label{eq:x-difference-expanded}
 \mathcal D_{\delta,\eta}f(x)
 ={}&\frac{(1+i\delta)(\eta-ix)}4
                [f(x+2i)-f(x)]\notag\\
 &+\frac{(1-i\delta)(\eta+ix)}4
                [f(x-2i)-f(x)].
\end{align}
Then
\begin{equation}\label{eq:x-eigenvalue}
 \mathcal D_{\delta,\eta}M_n=nM_n.
\end{equation}
The imaginary shifts are harmless on polynomials. To prove the
identity, factor the generating function as
\begin{equation}\label{eq:G-factorization}
 G(x;t)=A^{-a}B^{-b},\qquad
 A=1+(\delta+i)t,\quad B=1+(\delta-i)t,\quad
 a=\frac{\eta+ix}{2},\quad b=\frac{\eta-ix}{2}.
\end{equation}
Then $G(x+2i;t)/G(x;t)=A/B$ and
$G(x-2i;t)/G(x;t)=B/A$. Substitution in
\eqref{eq:x-difference-expanded} gives
\[
 \mathcal D_{\delta,\eta}G
 =-t\left(\frac{(\delta+i)a}{A}
          +\frac{(\delta-i)b}{B}\right)G
 =t\partial_tG.
\]
Coefficient comparison proves \eqref{eq:x-eigenvalue}.
Together with the three-term recurrence in $n$, this is the usual
bispectral structure of Meixner--Pollaczek polynomials
\cite[\S9.7]{KoekoekLeskySwarttouw2010}.

In the parameter direction, the lowering identity and differential
equation proved in Section~\ref{sec:parameter-duality} tell a different
story. The first-derivative coefficient in
\eqref{eq:delta-ODE-expanded} depends on $n$. This equation therefore
does not exhibit one fixed second-order operator in $\delta$ with
all $M_n$ as eigenfunctions. What we have established in the pair
$(x,\delta)$ is a precise parameter duality: an ordinary Jacobi
spectral problem in one direction and a positive
$R_{\mathrm{II}}$ pencil in the other. Calling this bispectrality
is useful provided that this distinction remains explicit.

\subsection{Why the Lie algebra has these coefficients}

There is a concise representation-theoretic reason for the factors
$n(n+\eta-1)$. On the basis $e_0,e_1,\ldots$ of
$\ell^2(\mathbb Z_{\ge0})$, define
\begin{align*}
 K_0e_n&=(n+\eta/2)e_n,\\
 K_+e_n&=\sqrt{(n+1)(n+\eta)}\,e_{n+1},\\
 K_-e_n&=\sqrt{n(n+\eta-1)}\,e_{n-1},\qquad K_-e_0=0.
\end{align*}
Direct calculation gives
\[
 [K_0,K_+]=K_+,\qquad [K_0,K_-]=-K_-,\qquad
 [K_+,K_-]=-2K_0,
\]
and $K_+^*=K_-$ on finite sequences. These are the positive
discrete-series relations for $\mathfrak{su}(1,1)$ with lowest
weight $\eta/2$. For arbitrary $\eta>0$, the group statement is
naturally made for the universal covering group; no integrality
condition is needed for the Lie algebra construction
\cite{hetyei2010meixner,koelink1998convolutions}.

The Jacobi operator is simply
\begin{equation}\label{eq:Lie-Jacobi}
 J_\delta=2\delta K_0+\sqrt{1+\delta^2}(K_++K_-).
\end{equation}
Normalize the polynomials by their norms:
\[
 p_n(x)=\frac{M_n(x;\delta,\eta)}
 {(1+\delta^2)^{n/2}\sqrt{n!(\eta)_n}}.
\]
The recurrence says exactly that the map $e_n\mapsto p_n$ takes
$J_\delta$ to multiplication by $x$ in
$L^2(\mu_\delta)$. This map is onto by the polynomial completeness proved in
Section~\ref{sec:meixner-time}.
The Jacobi Carleman condition
\[
 \sum_{n\ge1}\lambda_n^{-1/2}=\infty
\]
also gives essential self-adjointness on finite sequences
\cite{chihara1978introduction}. The density $K_\eta(x,\delta)$
is the spectral measure of the vector $e_0$. The numbers $p_n(x)$
are components of generalized eigenvectors; for fixed real $x$
they are not being asserted to form an $\ell^2$ eigenvector.

The parameter $\delta$ changes the operator inside the same
representation. Put $K_1=(K_++K_-)/2$,
$K_2=(K_+-K_-)/(2i)$, and $\delta=\sinh r$.
The commutators give
\[
 e^{-irK_2}K_1e^{irK_2}=\cosh r\,K_1+\sinh r\,K_0.
\]
One obtains this by differentiating the two conjugates of $K_1$
and $K_0$; they satisfy the two equations $X'=Y$, $Y'=X$.
The identity holds first on the common analytic core and then
for the self-adjoint closures. Hence
\begin{equation}\label{eq:Lie-conjugacy}
 J_\delta=e^{-irK_2}J_0e^{irK_2}.
\end{equation}
The operators are unitarily conjugate, while their spectral
measures at the fixed vector $e_0$ vary. This makes the exponential
tilt compatible with the unchanged continuous spectrum.

The wider element $2cK_0+K_++K_-$ leads, after elementary affine
changes, to Meixner-II, Laguerre, and Meixner-I according as
$|c|<1$, $|c|=1$, and $|c|>1$
\cite{hetyei2010meixner}. The same discriminant that appeared in
the Sheffer classification has reappeared in the Lie algebra.
Principal-series representations enter a different, larger story:
mixed tensor products of positive and negative discrete series,
continuous dual Hahn coefficients, and nonpolynomial Meixner
functions \cite{groenevelt2002meixner}. They are not needed to
explain the positive measure studied here.

\section{What carries over from moment duality?}\label{sec:outlook}

Ismail and Stanton's \emph{Moment duality for orthogonal polynomials}
develops parameter duality of type $R_{\mathrm I}$
\cite{IsmailStanton2025}. Their Meixner example varies the shape
parameter of the first-kind family. The coefficient multiplying
the preceding polynomial is linear in that parameter. Here the
conventional second-kind parametrization gives $1+\delta^2$,
so the corresponding problem is of type $R_{\mathrm{II}}$.
The analogy is real, but the analytic consequences are different.

In the $x$ direction, the density has exponential tails, and every
ordinary moment exists. In the parameter direction,
\begin{equation}\label{eq:parameter-tails}
 K_\eta(x,t)\sim C_\eta(x)e^{\pm\pi x/2}|t|^{-\eta}
       \qquad(t\to\pm\infty).
\end{equation}
For $\eta>1$, the ordinary absolute moments of $\nu_x$ exist
only in orders $k<\eta-1$; those of $\widehat\sigma_x$ exist
only in orders $k<\eta+1$. At fixed finite $\eta$, neither
probability provides all ordinary polynomial moments in the
parameter. What survives at every degree is the rational
orthogonality \eqref{eq:rational-orthogonality} or
\eqref{eq:weighted-rational-orthogonality}.

This observation also guides the combinatorics. The Motzkin model
in Section~\ref{sec:fractions} naturally describes an ordinary
Jacobi recurrence. Kim and Song's theory of generalized
$R_{\mathrm{II}}$ orthogonality provides the relevant broader
path model \cite{KimSong2026}: vertical-down and backward-down
steps join the usual up, level, and down steps, recording the
additional terms of the quadratic coefficient. Their formal
framework does not require a heavy-tailed positive measure to
have moments of all orders. Applying their model to the present
coefficients is a concrete combinatorial specialization, distinct
from the positive-measure statement proved here.

The new probabilistic readings sharpen this distinction. The finite
pencil already couples two exact Cauchy laws, before any limit is taken.
Its zero curves approach Laguerre rays, and their quadrature weights
approach the gamma law governing simultaneous extremes. At the other
end of the paper, independent increments explain why Carlitz's
continuous Hahn measure appears when a Meixner process is conditioned
on its endpoint. These are consequences of the same kernel, reached
by asking different questions of it.

There are still concrete questions to ask. What is the essential
spectrum of the Cauchy conditional operator $Q_\eta$? How quickly does
its Laguerre quadrature approximate the tail coefficient, especially
when the corner of $u\mapsto\min(1,u)$ lies close to a node? Which other
positive $R_{\mathrm{II}}$ pencils admit finite couplings with two
explicit marginals? The formulas above give starting points for these
questions.

The kernel gives a useful test for any further analogy. In one
direction it has exponential tails and ordinary orthogonality; in
the other it has algebraic tails and rational orthogonality. The
joint Cauchy law makes these two descriptions compatible without
erasing their difference. What looked like a mystery about the
number of measures is a question about how one kernel is read.

\appendix


\section{How much information survives conditioning?}\label{app:conditioning}

The two polynomial identities in
Proposition~\ref{prop:research-two-martingales} determine an entire
operator. Fix $\delta$ and $0<s<t$, retaining the notation
$\mu_t^{(\delta)}(dx)=K_t(x,\delta)\,dx$ from
Section~\ref{sec:meixner-time}. Define
\[
 T_{s,t}:L^2(\mu_t^{(\delta)})\longrightarrow L^2(\mu_s^{(\delta)}),
 \qquad
 (T_{s,t}f)(y)=\int_{\R}f(y+z)K_{t-s}(z,\delta)\,dz.
\]
This is conditional expectation from time $t$ to time $s$; its adjoint
is conditional expectation in the reverse direction. Both are
contractions by Jensen's inequality. The complete polynomial bases
make their action particularly simple.

\begin{proposition}\label{prop:research-conditional-spectrum}
Put $p_n^{(t)}=M_n^{(t)}/\sqrt{h_n(t)}$, with $h_n(t)$ as in
\eqref{eq:time-norms}. A complete singular system is
\begin{equation}\label{eq:research-singular-values}
 T_{s,t}p_n^{(t)}=\sigma_n p_n^{(s)},\qquad
 T_{s,t}^*p_n^{(s)}=\sigma_n p_n^{(t)},\qquad
 \boxed{\sigma_n=\sqrt{\frac{(s)_n}{(t)_n}}.}
\end{equation}
The operator is compact for every $t>s$, and, for $0<p<\infty$,
\begin{equation}\label{eq:research-schatten}
 T_{s,t}\in\mathcal S_p
 \quad\Longleftrightarrow\quad p(t-s)>2.
\end{equation}
Here $\mathcal S_p$ means that the $p$th powers of the singular values
are summable.
\end{proposition}

\begin{proof}
Normalize the two identities of
Proposition~\ref{prop:research-two-martingales}. The completeness
proved there gives the entire singular system. Moreover,
\[
 \frac{\sigma_{n+1}^2}{\sigma_n^2}=\frac{n+s}{n+t}<1,
 \qquad
 \sigma_n\sim\sqrt{\frac{\Gamma(t)}{\Gamma(s)}}\,
                    n^{-(t-s)/2}.
\]
The second assertion follows from the gamma ratio asymptotic. Thus the
singular values tend to zero, proving compactness; summing their
$p$th powers proves \eqref{eq:research-schatten}.
\end{proof}

In particular, the Hilbert--Schmidt threshold is $t-s>1$, and the
trace class threshold is $t-s>2$. On mean-zero functions the exact
norm is $\sigma_1=\sqrt{s/t}$. Equivalently, for square-integrable
$f,g$,
\begin{equation}\label{eq:time-maximal-correlation}
 |\operatorname{Cov}(f(X_s),g(X_t))|
 \le\sqrt{\frac{s}{t}}
       \sqrt{\operatorname{Var}(f(X_s))\operatorname{Var}(g(X_t))}.
\end{equation}
Equality is attained by the centered linear functions. The tilt changes
the laws and their orthogonal bases, but leaves every singular value
unchanged.

There is also an exact measure of dependence. For probability measures
$P\ll Q$, define
\[
 \chi^2(P\Vert Q)=\int\left(\frac{dP}{dQ}-1\right)^2\,dQ,
\]
allowing the value $+\infty$.

\begin{corollary}\label{cor:research-chi-square}
If $P_{s,t}$ is the joint law of $(X_s,X_t)$, then
\begin{equation}\label{eq:research-chi-square}
 \boxed{\displaystyle
 \chi^2(P_{s,t}\Vert\mu_s^{(\delta)}\otimes\mu_t^{(\delta)})
 =\begin{cases}
 s/(t-s-1),&t-s>1,\\
 +\infty,&0<t-s\le1.
 \end{cases}}
\end{equation}
When finite, $\|T_{s,t}\|_{\mathrm{HS}}^2=(t-1)/(t-s-1)$.
\end{corollary}

\begin{proof}
Relative to the product of the marginal laws, the joint density is
\[
 R(y,x)=\frac{K_{t-s}(x-y,\delta)}{K_t(x,\delta)}.
\]
For fixed $y$, the tail estimate \eqref{eq:Ceta-tail} gives
$R(y,x)=O(|x|^{-s})$ at infinity. Thus $R(y,\cdot)$ is bounded and
belongs to $L^2(\mu_t^{(\delta)})$. Its $n$th orthogonal coefficient
is $T_{s,t}p_n^{(t)}(y)=\sigma_n p_n^{(s)}(y)$.
Parseval followed by Tonelli gives, with infinite values allowed,
\[
 \int R^2\,d(\mu_s^{(\delta)}\otimes\mu_t^{(\delta)})
 =\sum_{n\ge0}\sigma_n^2.
\]
The constant term is one, leaving
$\chi^2=\sum_{n\ge1}\sigma_n^2$.

To sum the series, put $r=t-s$ and let $B$ have beta density
$u^{s-1}(1-u)^{r-1}/\mathrm B(s,r)$ on $(0,1)$. Then
\[
 \sum_{n\ge0}\frac{(s)_n}{(t)_n}
 =\E\sum_{n\ge0}B^n
 =\E\frac1{1-B}
 =\frac{\mathrm B(s,r-1)}{\mathrm B(s,r)}
 =\frac{t-1}{r-1}
\]
if $r>1$. If $r\le1$, the same integral diverges at $1$.
This proves both assertions.
\end{proof}

The associated density expansion is an $L^2$ expansion precisely when
$t-s>1$; the operator identities \eqref{eq:research-singular-values}
hold for every positive time gap. Thus the amount of smoothing is
measured exactly by elapsed time, degree by degree.


\section{What the zero curves approach}\label{sec:zero-geometry}

The positive pencil does more than keep the zeros real: it identifies
their directions at infinity. The eigenvalues of $A_N$ are the zeros
of $L_N^{(\eta-1)}$. Thus the Laguerre family appears as the asymptotic
geometry of the same curves that support the finite coupling in
Theorem~\ref{thm:finite-coupling}.

\subsection{A first correction to the Laguerre limit}

Fix $N\ge1$ and $\eta>0$, and write
\[
 \mathcal L_N(y)=(-1)^N N!L_N^{(\eta-1)}(y),\qquad
 0<\ell_1<\cdots<\ell_N
\]
for the monic Laguerre polynomial and its zeros. The leading
Meixner--Pollaczek-to-Laguerre limit is classical; the corresponding
reciprocal-Laguerre approximation for the parameter zeros appears in
\cite[Theorem~4.2]{SilvaRibeiroEtAl2025}. Full polynomial expansions
are available in \cite[\S4]{TemmeLopez2001}. In the present normalization,
the first correction has a particularly short derivation and an explicit
expression involving only the limiting zero.

\begin{proposition}\label{prop:laguerre-correction}
Let $x_1(\delta)<\cdots<x_N(\delta)$ be the zeros of
$M_N(x;\delta,\eta)$ in $x$. As $\delta\to+\infty$,
\begin{equation}\label{eq:zero-laguerre-correction}
 x_j(\delta)=\delta\ell_j+\frac{\beta_j}{\delta}
                 +O(\delta^{-3}),\qquad
 \beta_j=\frac{2\ell_j-2N-\eta+\eta(2-\eta)/\ell_j}{6}.
\end{equation}
If $d_1(x)<\cdots<d_N(x)$ are the zeros in the parameter, then
\begin{equation}\label{eq:parameter-laguerre-correction}
 d_k(x)=\frac{x}{\ell_{N+1-k}}
             -\frac{\beta_{N+1-k}}{x}+O(x^{-3}),
             \qquad x\to+\infty.
\end{equation}
All error estimates here are for fixed $N$ and $\eta$.
\end{proposition}

\begin{proof}
Put $h=\delta^{-2}$. Rescale the generating function by replacing
$x$ with $\delta y$ and its series variable with $u/\delta$. This gives
\begin{align*}
 \sum_{n\ge0}\delta^{-n}M_n(\delta y;\delta,\eta)\frac{u^n}{n!}
 &=[(1+u)^2+hu^2]^{-\eta/2}
    \exp\!\left\{\frac{y}{\sqrt h}
            \arctan\frac{\sqrt h\,u}{1+u}\right\}\\
 &=F_0(y,u)\left[1-h\left(\frac{\eta r^2}{2}
                         +\frac{yr^3}{3}\right)+O(h^2)\right],
\end{align*}
where $r=u/(1+u)$ and
\[
 F_0(y,u)=(1+u)^{-\eta}e^{yr}
       =\sum_{n\ge0}\mathcal L_n(y)\frac{u^n}{n!}.
\]
The apparent square root causes no difficulty: the displayed expression
is analytic in $h$ near zero. Since $\partial_y^mF_0=(u/(1+u))^mF_0$
for every nonnegative integer $m$, coefficient comparison yields
\begin{equation}\label{eq:polynomial-laguerre-correction}
 \delta^{-N}M_N(\delta y;\delta,\eta)
 =\mathcal L_N(y)-\delta^{-2}
    \left(\frac\eta2\mathcal L_N''(y)
              +\frac y3\mathcal L_N'''(y)\right)+O(\delta^{-4}).
\end{equation}
This is a coefficientwise polynomial expansion. Each simple Laguerre
zero therefore moves analytically with $h$, and its first displacement is
\[
 \beta_j=\frac\eta2\frac{\mathcal L_N''(\ell_j)}{\mathcal L_N'(\ell_j)}
       +\frac{\ell_j}{3}
             \frac{\mathcal L_N'''(\ell_j)}{\mathcal L_N'(\ell_j)}.
\]
The Laguerre equation
$y\mathcal L_N''+(\eta-y)\mathcal L_N'+N\mathcal L_N=0$
and its derivative give, at $y=\ell_j$,
\[
 \frac{\mathcal L_N''}{\mathcal L_N'}=\frac{\ell_j-\eta}{\ell_j},
 \qquad
 \frac{\mathcal L_N'''}{\mathcal L_N'}
 =\frac{(\ell_j-\eta-1)(\ell_j-\eta)-\ell_j(N-1)}{\ell_j^2}.
\]
Substitution simplifies to \eqref{eq:zero-laguerre-correction}.
Finally, inverting $x=\ell_j\delta+\beta_j/\delta+O(\delta^{-3})$
gives $\delta=x/\ell_j-\beta_j/x+O(x^{-3})$.
The order reverses because the positive numbers $1/\ell_j$ decrease.
\end{proof}

\begin{remark}
Reflection, $M_N(-x;-\delta,\eta)=(-1)^NM_N(x;\delta,\eta)$,
gives the negative-infinity formulas with the indices reversed.
\end{remark}

\subsection{The angle reveals a convex equilibrium}

The parameter zeros also admit an electrostatic interpretation, already
established for complementary Romanovski--Routh polynomials in
\cite[\S3]{MartinezFinkelshteinEtAl2019}. The substitution $d_j=\tan\theta_j$
makes both uniqueness and motion of the equilibrium particularly clear.

\begin{proposition}\label{prop:angular-equilibrium}
The angles $\theta_j(x)=\arctan d_j(x)$ uniquely minimize
\begin{equation}\label{eq:angular-energy}
 E_x(\theta)=-2\sum_{j<k}\log\sin(\theta_k-\theta_j)
            -\eta\sum_j\log\cos\theta_j-x\sum_j\theta_j
\end{equation}
on $-\pi/2<\theta_1<\cdots<\theta_N<\pi/2$. Moreover,
\begin{equation}\label{eq:angular-speed}
 0<\theta_j'(x)=\frac{d_j'(x)}{1+d_j(x)^2}\le\frac1\eta.
\end{equation}
The upper bound is strict when $N>1$.
\end{proposition}

\begin{proof}
At a parameter zero, the differential equation
\eqref{eq:delta-ODE-expanded} gives
\[
 2\sum_{k\ne j}\frac1{d_j-d_k}
       =\frac{(\eta+2N-2)d_j-x}{1+d_j^2}.
\]
Using $\cot(\theta_j-\theta_k)
=(1+d_j^2)/(d_j-d_k)-d_j$, this becomes
\[
 2\sum_{k\ne j}\cot(\theta_j-\theta_k)
       =\eta\tan\theta_j-x,
\]
which is exactly $\nabla E_x=0$. If $H$ is the Hessian, then
\[
 v^{\mathsf T}Hv
 =\eta\sum_jv_j^2\sec^2\theta_j
       +2\sum_{j<k}(v_j-v_k)^2\csc^2(\theta_j-\theta_k)>0
       \qquad(v\ne0).
\]
Thus $E_x$ is strictly convex. It tends to infinity at every boundary
of the ordered chamber, proving existence and uniqueness.

Differentiating the equilibrium equations gives $H\theta'=\mathbf1$.
Write $u=\theta'$. At an index where $u_j$ is smallest,
\[
 1=(Hu)_j
   =\eta\sec^2\theta_j\,u_j
       +2\sum_{k\ne j}\csc^2(\theta_j-\theta_k)(u_j-u_k)
   \le\eta\sec^2\theta_j\,u_j,
\]
so every component is positive. At a largest component the inequality
reverses, giving $u_j\le\cos^2\theta_j/\eta\le1/\eta$.
Equality in the last bound would force all components of $u$ to agree
and all angles to vanish, which is impossible for $N>1$.
\end{proof}

The energy describes repelling points confined between two endpoints.
Increasing $x$ applies a uniform force to the right; strict convexity
ensures a unique response, and \eqref{eq:angular-speed} bounds the speed
of every point.

\section*{Acknowledgement} The author is grateful to Maxim Derevyagin for asking about bispectrality of Meixner-II polynomials and for pointing out his relevant work. 

\section*{AI Declaration} Artificial intelligence (AI) tools were used to generate the graphics for this paper. An LLM was also used to audit its earlier draft and assist in the rewriting. The author assumes full responsibility for all results presented here.

\bigskip
\noindent{\small Department of Mathematics, University of California, Berkeley}


\begin{thebibliography}{41}
\providecommand{\natexlab}[1]{#1}
\providecommand{\url}[1]{\texttt{#1}}
\expandafter\ifx\csname urlstyle\endcsname\relax
  \providecommand{\doi}[1]{doi: #1}\else
  \providecommand{\doi}{doi: \begingroup \urlstyle{rm}\Url}\fi

\bibitem[Bateman(1933)]{bateman1933properties}
Harry Bateman.
\newblock Some properties of a certain set of polynomials.
\newblock \emph{Tohoku Mathematical Journal, First Series}, 37:\penalty0
  23--38, 1933.
\newblock URL
  \url{https://www.jstage.jst.go.jp/article/tmj1911/37/0/37_0_23/_article/-char/en}.

\bibitem[Bracciali et~al.(2020)Bracciali, Pereira, and
  Sri~Ranga]{BraccialiEtAl2020}
C.~F. Bracciali, J.~A. Pereira, and A.~Sri~Ranga.
\newblock Quadrature rules from a {$R_{II}$} type recurrence relation and
  associated quadrature rules on the unit circle.
\newblock \emph{Numerical Algorithms}, 83:\penalty0 1029--1061, 2020.
\newblock \doi{10.1007/s11075-019-00714-w}.

\bibitem[Brent and Harvey(2013)]{BrentHarvey2013}
Richard~P. Brent and David Harvey.
\newblock Fast computation of {B}ernoulli, tangent and secant numbers.
\newblock In \emph{Computational and analytical mathematics}, volume~50 of
  \emph{Springer Proc. Math. Stat.}, pages 127--142. Springer, New York, 2013.
\newblock ISBN 978-1-4614-7621-4; 978-1-4614-7620-7.
\newblock \doi{10.1007/978-1-4614-7621-4_8}.

\bibitem[Butzer and Koornwinder(2019)]{butzer2019josef}
Paul~L. Butzer and Tom~H. Koornwinder.
\newblock Josef {M}eixner: his life and his orthogonal polynomials.
\newblock \emph{Indag. Math. (N.S.)}, 30\penalty0 (1):\penalty0 250--264, 2019.
\newblock ISSN 0019-3577,1872-6100.
\newblock \doi{10.1016/j.indag.2018.09.009}.
\newblock URL \url{https://arxiv.org/abs/1609.02588}.
\newblock Corrected version: arXiv:1609.02588v5 (4 January 2024); see also the
  authors' errata of 5 January 2024.

\bibitem[Carlitz(1959)]{carlitz1959bernoulli}
L.~Carlitz.
\newblock Bernoulli and {E}uler numbers and orthogonal polynomials.
\newblock \emph{Duke Mathematical Journal}, 26\penalty0 (1):\penalty0 1--15,
  1959.
\newblock \doi{10.1215/S0012-7094-59-02601-8}.

\bibitem[Chihara(1978)]{chihara1978introduction}
Theodore~S. Chihara.
\newblock \emph{An Introduction to Orthogonal Polynomials}.
\newblock Gordon and Breach Science Publishers, New York/London/Paris, 1978.

\bibitem[Derevyagin(2017)]{Derevyagin2017}
Maxim Derevyagin.
\newblock A note on {Wall}'s modification of the {Schur} algorithm and linear
  pencils of {Jacobi} matrices.
\newblock \emph{Journal of Approximation Theory}, 221:\penalty0 1--21, 2017.
\newblock \doi{10.1016/j.jat.2017.05.001}.

\bibitem[Devroye and Hill(2025)]{DevroyeHill2025}
Luc Devroye and J.~R. Hill.
\newblock Simulating random variates from the {Pearson IV} and betaized
  {Meixner--Morris} distributions, 2025.
\newblock URL \url{https://arxiv.org/abs/2505.15930}.
\newblock Preprint, arXiv:2505.15930.

\bibitem[Flajolet(1980)]{Flajolet1980}
Philippe Flajolet.
\newblock Combinatorial aspects of continued fractions.
\newblock \emph{Discrete Mathematics}, 32\penalty0 (2):\penalty0 125--161,
  1980.
\newblock \doi{10.1016/0012-365X(80)90050-3}.

\bibitem[Groenevelt and Koelink(2002)]{groenevelt2002meixner}
Wolter Groenevelt and Erik Koelink.
\newblock Meixner functions and polynomials related to {L}ie algebra
  representations.
\newblock \emph{Journal of Physics A: Mathematical and General}, 35\penalty0
  (1):\penalty0 65--85, 2002.

\bibitem[Gupta and Milovanovi{\'c}(2022)]{GuptaMilovanovic2022}
Vijay Gupta and Gradimir~V. Milovanovi{\'c}.
\newblock A solution to exponential operators.
\newblock \emph{Results in Mathematics}, 77, 2022.
\newblock \doi{10.1007/s00025-022-01743-z}.
\newblock Article 207.

\bibitem[Hetyei(2010)]{hetyei2010meixner}
G\'abor Hetyei.
\newblock Meixner polynomials of the second kind and quantum algebras
  representing {${\rm su}(1,1)$}.
\newblock \emph{Proc. R. Soc. Lond. Ser. A Math. Phys. Eng. Sci.}, 466\penalty0
  (2117):\penalty0 1409--1428, 2010.
\newblock ISSN 1364-5021,1471-2946.
\newblock \doi{10.1098/rspa.2009.0497}.

\bibitem[Holtz(2024)]{Holtz2024Zeta}
Olga Holtz.
\newblock Classical values of {Z}eta, as simple as possible but not simpler.
\newblock \emph{Notices of the American Mathematical Society}, 71\penalty0
  (6):\penalty0 725--731, 2024.
\newblock \doi{10.1090/noti2949}.

\bibitem[Ismail and Masson(1995)]{IsmailMasson1995}
Mourad E.~H. Ismail and D.~R. Masson.
\newblock Generalized orthogonality and continued fractions.
\newblock \emph{Journal of Approximation Theory}, 83\penalty0 (1):\penalty0
  1--40, 1995.
\newblock \doi{10.1006/jath.1995.1106}.

\bibitem[Ismail and May(1978)]{IsmailMay1978}
Mourad E.~H. Ismail and C.~P. May.
\newblock On a family of approximation operators.
\newblock \emph{Journal of Mathematical Analysis and Applications},
  63:\penalty0 446--462, 1978.
\newblock \doi{10.1016/0022-247X(78)90090-2}.

\bibitem[Ismail and Sri~Ranga(2019)]{IsmailSriRanga2019}
Mourad E.~H. Ismail and A.~Sri~Ranga.
\newblock {$R_{II}$} type recurrence, generalized eigenvalue problem and
  orthogonal polynomials on the unit circle.
\newblock \emph{Linear Algebra and its Applications}, 562:\penalty0 63--90,
  2019.
\newblock \doi{10.1016/j.laa.2018.10.005}.

\bibitem[Ismail and Stanton(2025)]{IsmailStanton2025}
Mourad E.~H. Ismail and Dennis Stanton.
\newblock Moment duality for orthogonal polynomials.
\newblock \emph{The Ramanujan Journal}, 68:\penalty0 54, 2025.
\newblock \doi{10.1007/s11139-025-01203-6}.
\newblock Article 54.

\bibitem[Kim and Song(2026)]{KimSong2026}
Jang~Soo Kim and Minho Song.
\newblock Combinatorics of generalized orthogonal polynomials of type
  {$R_{II}$}.
\newblock \emph{The Electronic Journal of Combinatorics}, 33\penalty0
  (4):\penalty0 P4.24, 2026.
\newblock \doi{10.37236/13718}.
\newblock Paper No. P4.24.

\bibitem[Koekoek et~al.(2010)Koekoek, Lesky, and
  Swarttouw]{KoekoekLeskySwarttouw2010}
Roelof Koekoek, Peter~A. Lesky, and Ren{\'e}~F. Swarttouw.
\newblock \emph{Hypergeometric Orthogonal Polynomials and Their
  {$q$}-Analogues}.
\newblock Springer Monographs in Mathematics. Springer, Berlin, Heidelberg,
  2010.
\newblock ISBN 978-3-642-05013-8.
\newblock \doi{10.1007/978-3-642-05014-5}.

\bibitem[Koelink and Van~der Jeugt(1998)]{koelink1998convolutions}
Erik Koelink and Joris Van~der Jeugt.
\newblock Convolutions for orthogonal polynomials from {L}ie and quantum
  algebra representations.
\newblock \emph{SIAM Journal on Mathematical Analysis}, 29\penalty0
  (3):\penalty0 794--822, 1998.

\bibitem[Koelink(1996)]{Koelink1996}
H.~T. Koelink.
\newblock On {J}acobi and continuous {H}ahn polynomials.
\newblock \emph{Proceedings of the American Mathematical Society}, 124\penalty0
  (3):\penalty0 887--898, 1996.
\newblock \doi{10.1090/S0002-9939-96-03190-5}.

\bibitem[Liaw and Treil(2020)]{LiawTreil2020}
Constanze Liaw and Sergei Treil.
\newblock Matrix measures and finite rank perturbations of self-adjoint
  operators.
\newblock \emph{Journal of Spectral Theory}, 10:\penalty0 1173--1210, 2020.
\newblock \doi{10.4171/JST/324}.
\newblock URL \url{https://arxiv.org/abs/1806.08856}.

\bibitem[Mart{\'i}nez-Finkelshtein et~al.(2019)Mart{\'i}nez-Finkelshtein,
  Silva~Ribeiro, Sri~Ranga, and Tyaglov]{MartinezFinkelshteinEtAl2019}
Andrei Mart{\'i}nez-Finkelshtein, L.~L. Silva~Ribeiro, A.~Sri~Ranga, and
  M.~Tyaglov.
\newblock Complementary {Romanovski--Routh} polynomials: From orthogonal
  polynomials on the unit circle to {Coulomb} wave functions.
\newblock \emph{Proceedings of the American Mathematical Society}, 147\penalty0
  (6):\penalty0 2625--2640, 2019.
\newblock \doi{10.1090/proc/14423}.

\bibitem[Mazzola and Muliere(2011)]{MazzolaMuliere2011}
E.~Mazzola and P.~Muliere.
\newblock Reviewing alternative characterizations of {Meixner} process.
\newblock \emph{Probability Surveys}, 8:\penalty0 127--154, 2011.
\newblock \doi{10.1214/11-PS177}.

\bibitem[Meixner(1934)]{meixner1934orthogonale}
J.~Meixner.
\newblock Orthogonale {P}olynomsysteme mit einer besonderen {G}estalt der
  erzeugenden {F}unktion.
\newblock \emph{Journal of the London Mathematical Society}, 9\penalty0
  (1):\penalty0 6--13, 1934.
\newblock \doi{10.1112/jlms/s1-9.1.6}.

\bibitem[Meixner(1933)]{Meixner1933Green}
Josef Meixner.
\newblock Die {G}reensche {F}unktion des wellenmechanischen {K}eplerproblems.
\newblock \emph{Mathematische Zeitschrift}, 36:\penalty0 677--707, 1933.
\newblock \doi{10.1007/BF01188644}.
\newblock Munich doctoral dissertation, 1931.

\bibitem[Meixner and Reik(1959)]{MeixnerReik1959}
Josef Meixner and H.~G. Reik.
\newblock Thermodynamik der irreversiblen {P}rozesse.
\newblock In Siegfried Fl{\"u}gge, editor, \emph{Prinzipien der {T}hermodynamik
  und {S}tatistik}, volume III/2 of \emph{Handbuch der Physik}, pages 413--523.
  Springer, Berlin, G{\"o}ttingen, Heidelberg, 1959.
\newblock \doi{10.1007/978-3-642-45912-2_4}.

\bibitem[Meixner and Sch{\"a}fke(1954)]{MeixnerSchaefke1954}
Josef Meixner and Friedrich~Wilhelm Sch{\"a}fke.
\newblock \emph{Mathieusche {F}unktionen und {S}ph{\"a}roidfunktionen mit
  {A}nwendungen auf physikalische und technische {P}robleme}, volume~71 of
  \emph{Die Grundlehren der mathematischen Wissenschaften}.
\newblock Springer, Berlin, G{\"o}ttingen, Heidelberg, 1954.
\newblock \doi{10.1007/978-3-662-00941-3}.

\bibitem[{National Institute of Standards and Technology}(n.d.)]{DLMF}
{National Institute of Standards and Technology}.
\newblock {NIST Digital Library of Mathematical Functions}.
\newblock \url{https://dlmf.nist.gov/18}, n.d.
\newblock URL \url{https://dlmf.nist.gov/}.
\newblock Chapters 5, 15, and 18: Gamma Function, Hypergeometric Function, and
  Orthogonal Polynomials.

\bibitem[Pollaczek(1950)]{Pollaczek1950Limits}
F{\'e}lix Pollaczek.
\newblock Sur une famille de polyn{\^o}mes orthogonaux qui contient les
  polyn{\^o}mes d'{H}ermite et de {L}aguerre comme cas limites.
\newblock \emph{Comptes Rendus de l'Acad{\'e}mie des Sciences de Paris},
  230:\penalty0 1563--1565, 1950.
\newblock URL
  \url{https://webusers.imj-prg.fr/~leila.schneps/grothendieckcircle/AG/AG-2.pdf}.

\bibitem[Schoutens and Teugels(1998)]{SchoutensTeugels1998}
Wim Schoutens and Jozef~L. Teugels.
\newblock L{\'e}vy processes, polynomials and martingales.
\newblock \emph{Communications in Statistics---Stochastic Models}, 14:\penalty0
  335--349, 1998.
\newblock \doi{10.1080/15326349808807475}.

\bibitem[Sheffer(1939)]{Sheffer1939}
I.~M. Sheffer.
\newblock Some properties of polynomial sets of type zero.
\newblock \emph{Duke Mathematical Journal}, 5\penalty0 (3):\penalty0 590--622,
  1939.
\newblock \doi{10.1215/S0012-7094-39-00549-1}.

\bibitem[Silva~Ribeiro et~al.(2025)Silva~Ribeiro, Sri~Ranga, and {Yen Chi
  Lun}]{SilvaRibeiroEtAl2025}
L.~L. Silva~Ribeiro, A.~Sri~Ranga, and {Yen Chi Lun}.
\newblock Complementary {Romanovski--Routh} polynomials and their zeros.
\newblock \emph{Trends in Computational and Applied Mathematics}, 26, 2025.
\newblock \doi{10.5540/tcam.2025.026.e01866}.
\newblock URL \url{https://arxiv.org/abs/2212.02260}.
\newblock Article e01866.

\bibitem[Sommerfeld(1952)]{Sommerfeld1952}
Arnold Sommerfeld.
\newblock \emph{Thermodynamik und {S}tatistik}, volume~5 of \emph{Vorlesungen
  {\"u}ber theoretische Physik}.
\newblock Dieterich, Wiesbaden, 1952.
\newblock Edited by Fritz Bopp and Josef Meixner.

\bibitem[Temme and L{\'o}pez(2001)]{TemmeLopez2001}
Nico~M. Temme and Jos{\'e}~L. L{\'o}pez.
\newblock The {Askey} scheme for hypergeometric orthogonal polynomials viewed
  from asymptotic analysis.
\newblock \emph{Journal of Computational and Applied Mathematics},
  133:\penalty0 623--633, 2001.
\newblock \doi{10.1016/S0377-0427(00)00683-X}.

\bibitem[Tojo and Yoshino(2026)]{TojoYoshino2026}
K.~Tojo and T.~Yoshino.
\newblock A family of conjugate priors for the natural exponential family
  generated by the generalized hyperbolic secant distributions.
\newblock \emph{Statistics \& Probability Letters}, 235:\penalty0 110739, 2026.
\newblock \doi{10.1016/j.spl.2026.110739}.

\bibitem[Truksa(1931)]{truksa1931}
L.~Truksa.
\newblock Hypergeometric orthogonal systems of polynomials {I}, {II}, {III}.
\newblock \emph{Aktu{\'a}rsk{\'e} v{\v{e}}dy}, 2:\penalty0 65--84, 113--144,
  177--203, 1931.
\newblock URL \url{http://dml.cz/dmlcz/144552}.

\bibitem[Viennot(1983)]{viennot1983combinatorial}
G{\'e}rard Viennot.
\newblock \emph{Une th{\'e}orie combinatoire des polyn{\^o}mes orthogonaux
  g{\'e}n{\'e}raux}.
\newblock Universit{\'e} du Qu{\'e}bec {\`a} Montr{\'e}al, Montr{\'e}al, 1983.
\newblock URL
  \url{https://www.imsc.res.in/~viennot/wa_files/viennotop1983-ocr.pdf}.
\newblock Lecture notes, September--October 1983.

\bibitem[Wall(1944)]{Wall1944}
H.~S. Wall.
\newblock Continued fractions and bounded analytic functions.
\newblock \emph{Bulletin of the American Mathematical Society}, 50:\penalty0
  110--119, 1944.
\newblock \doi{10.1090/S0002-9904-1944-08094-4}.

\bibitem[Wall(1946)]{Wall1946}
H.~S. Wall.
\newblock Bounded {$J$}-fractions.
\newblock \emph{Bulletin of the American Mathematical Society}, 52:\penalty0
  686--693, 1946.

\bibitem[Wall(1948)]{Wall1948}
H.~S. Wall.
\newblock \emph{Analytic Theory of Continued Fractions}.
\newblock D. Van Nostrand Company, New York, 1948.

\end{thebibliography}
\end{document}